\documentclass[11pt]{article}
\usepackage[centertags]{amsmath}
\usepackage{pdfsync,bbm,color} 
\hsize \textwidth
\advance \hsize by -\marginparwidth
\evensidemargin
\oddsidemargin
\usepackage{amssymb}
\usepackage{graphicx}

 \usepackage[active]{srcltx}

\usepackage[pagebackref,colorlinks=true,pdfpagemode=none,urlcolor=blue,linkcolor=blue,citecolor=blue]{hyperref}

\advance\hoffset by 5mm

\usepackage{amssymb}

\newtheorem{df}{Definition}[section]
\newtheorem{lemma}[df]{Lemma}
\newtheorem{prop}[df]{Proposition}

\newtheorem{thm}[df]{Theorem}

\newtheorem{e}[df]{Example}
\newtheorem{cor}[df]{Corollary}

\newcommand{\pdr}[2]{\dfrac{\partial{#1}}{\partial{#2}}}

\newcommand{\pdrr}[2]{\dfrac{\partial^2{#1}}{\partial{#2}^2}}

\makeatletter \@addtoreset{equation}{section}

\newcommand{\eps}{\varepsilon}
\newcommand{\bes}{\begin{displaymath}}
\newcommand{\ees}{\end{displaymath}}
\newcommand{\be}{\begin{equation}}
\newcommand{\ee}{\end{equation}}
\newcommand{\ba}{\begin{eqnarray}}
\newcommand{\ea}{\end{eqnarray}}
\newcommand{\bas}{\begin{eqnarray*}}
\newcommand{\eas}{\end{eqnarray*}}

\newcommand{\bal}{\begin{aligned}}
\newcommand{\enbal}{\end{aligned}}

\newcommand{\@Bbb}[1]{\ensuremath{\mathbb #1}}

\newcommand{\B}{{\@Bbb B}}
\newcommand{\C}{{\@Bbb C}}
\newcommand{\E}{{\@Bbb E}}
\newcommand{\F}{{\@Bbb F}}
\renewcommand{\P}{{\@Bbb P}}

\newcommand{\Q}{{\@Bbb Q}}
\newcommand{\bQ}{{\@Bbb Q}}
\newcommand{\N}{{\@Bbb N}}
\newcommand{\R}{{\@Bbb R}}
\newcommand{\bbR}{{\@Bbb R}}
\newcommand{\W}{{\@Bbb W}}
\newcommand{\Z}{{\@Bbb Z}}
\newcommand{\bbZ}{{\@Bbb Z}}

\newcommand{\calC}{{\mathcal C}}

\newcommand{\Rm}{{\mathbb R}}
\newcommand{\Nm}{{\mathbb N}}

\newcommand{\Em}{{\mathbb E}}

\newcommand{\@s}[1]{\ensuremath{\mathcal #1}}
\newcommand{\cA}{\@s A}
\newcommand{\cB}{\@s B}
\newcommand{\cC}{\@s C}
\newcommand{\cD}{\@s D}
\newcommand{\cE}{\@s E}
\newcommand{\cF}{\@s F}
\newcommand{\cG}{\@s G}
\newcommand{\cH}{\@s H}
\newcommand{\cI}{\@s I}
\newcommand{\cJ}{\@s J}
\newcommand{\cK}{\@s K}
\newcommand{\cL}{\@s L}
\newcommand{\cN}{\@s N}
\newcommand{\cM}{\@s M}
\newcommand{\cO}{\@s O}
\newcommand{\cP}{\@s P}
\newcommand{\cR}{\@s R}
\newcommand{\cS}{\@s S}
\newcommand{\cT}{\@s T}
\newcommand{\cV}{\@s V}
\newcommand{\cW}{\@s W}
\newcommand{\cX}{\@s X}
\newcommand{\cY}{\@s Y}
\newcommand{\cZ}{\@s Z}

\newcommand{\@bm}[1]{\ensuremath{\mathbf #1}}
\newcommand{\bma}{\@bm a}
\newcommand{\bmb}{\@bm b}
\newcommand{\bmc}{\@bm c}
\newcommand{\bmd}{\@bm d}

\newcommand{\bme}{\@bm e}
\newcommand{\bmf}{\@bm f}
\newcommand{\bmg}{\@bm g}
\newcommand{\bmh}{\@bm h}
\newcommand{\bmi}{\@bm i}
\newcommand{\bmj}{\@bm j}
\newcommand{\bmk}{\@bm k}
\newcommand{\bml}{\@bm l}
\newcommand{\bmm}{\@bm m}
\newcommand{\bmn}{\@bm n}
\newcommand{\bmo}{\@bm o}
\newcommand{\bmp}{\@bm p}
\newcommand{\bmq}{\@bm q}
\newcommand{\bmr}{\@bm r}
\newcommand{\bms}{\@bm s}
\newcommand{\bmt}{\@bm t}
\newcommand{\bmu}{\@bm u}
\newcommand{\bmw}{\@bm w}
\newcommand{\bmv}{\@bm v}
\newcommand{\bmx}{\@bm x}
\newcommand{\bx}{\@bm x}
\newcommand{\bmy}{\@bm y}
\newcommand{\bmz}{\@bm z}

\newcommand{\by}{\@bm y}
\newcommand{\bmzero}{\@bm 0}

\newcommand{\@g}[1]{\ensuremath{\mathfrak #1}}
\newcommand{\gA}{\@g A}
\newcommand{\gD}{\@g D}
\newcommand{\gJ}{\@g J}
\newcommand{\gF}{\@g F}
\newcommand{\gM}{\@g M}
\newcommand{\gR}{\@g R}

\newcommand{\commentout}[1]{{}}

\newcommand{\farc}{\frac}

\definecolor{violet}{rgb}{0.58, 0.0, 0.83}

\newcommand{\Pm}{\mathbb{P}}

\def\di{\displaystyle}

\definecolor{Red}{rgb}{1,0,0}

\newcommand{\bbL}{\mathbb L}
\newcommand{\qed}{$\Box$}

\newcommand{\one}{\mathbbm{1}}

\makeatother
\begin{document}

\author{Leonid Mytnik\footnote{Faculty of Industrial Engineering and Management, Technion, Technion City, Haifa 3200003, Israel; email: leonid@ie.technion.ac.il}
\and Jean-Michel Roquejoffre\footnote{Institut de Math\'ematiques de Toulouse; UMR 5219, Universit\'e de Toulouse; CNRS, UPS IMT, F-31062 Toulouse Cedex 9, France; 
e-mail: Jean-Michel.Roquejoffre@math.univ-toulouse.fr}
\and Lenya Ryzhik\footnote{Department of Mathematics, Stanford University, Stanford CA 94350, USA; email: ryzhik@stanford.edu}
}
\title{
Fisher-KPP equation with small data and the extremal process of branching Brownian motion}
 
\maketitle

\begin{abstract}
We consider the limiting 
extremal process $\cX$ of the particles of the binary branching Brownian motion. 
We show that after a shift by the logarithm of the derivative martingale $Z$, the rescaled "density" of particles, which are
at  distance $n+x$ from a position close to the tip of $\cX$,
 converges in probability  to a multiple of the exponential $e^x$  as $n\to+\infty$. 
We also show
that the fluctuations of the density, after another scaling and an additional random but explicit shift, 
converge to a~$1$-stable random variable. Our approach uses analytic techniques and is motivated by the connection  between
the properties of the branching Brownian motion and the Bramson shift of the  solutions 
to the Fisher-KPP equation with some specific
initial conditions initiated in~\cite{BD1,BD2} and further developed in the present paper. 
The proofs of the limit theorems for~$\cX$    
rely crucially on the fine asymptotics of the  behavior of the 
Bramson shift for  the Fisher-KPP equation starting with initial conditions of "size" 
$0<\eps\ll 1$,  up to terms of the order $[{(\log \eps^{-1})]^{-1-\gamma}}$, with some~$\gamma>0$. 

\end{abstract}
\section{Introduction}

\subsection*{The BBM connection to the Fisher-KPP equation}

The standard binary branching Brownian motion  (BBM) on $\Rm$ is a process that starts at the initial
time $t=0$ with one particle at the position $x=0$. The particle performs a Brownian motion until a random,
exponentially distributed with parameter $1$, 
time $\tau>0$ when it splits into two off-spring particles. 
Each of the  new particles 
starts an independent Brownian motion at the branching point. They have independent, 
again exponentially distributed with parameter $1$, clocks attached
to them that ring at their respective branching times. At the branching time,
the particle
that is branching produces two independent off-spring,
and the process continues, so that at a time~$t>0$ we have a random number $N_t$ of Brownian particles. It will be convenient for us to assume that the variance of each individual Brownian motion
is not $t$ but~$2t$. 

A remarkable observation by McKean \cite{McK} is a connection between the location 
of the maximal particle for the BBM and the Fisher-KPP equation
\begin{equation}\label{20apr1402}
\pdr{u}t=\frac{\partial^2 u}{\partial x^2}+u -u^2,
\end{equation}
with the initial condition $u(0,x)=\mathbbm{1}(x<0)$. 
Let~$x_1(t)\geq x_2(t)\geq \ldots\ge x_{N_t}(t)$ be the positions of the
BBM particles at time $t$. McKean has discovered an exact formula
\be\label{20apr1404}
u(t,x)=\Pm[x_1(t)>x].
\ee
More generally, given a sufficiently regular function $g(x)$, the 
solution to the Fisher-KPP equation~(\ref{20apr1402}) with the initial condition $u(0,x)=g(x)$ can be written as  
\be\label{20apr1406}
\bal
u(t,x)&=1-\Em_x\Big[\prod_{j=1}^{N_t}(1-g(x_j(t))\Big]
=1-\Em_0\Big[\prod_{j=1}^{N_t}(1-g(x+x_j(t))\Big]\\
&=1-\Em_0\Big[\prod_{j=1}^{N_t}(1-g(x-x_j(t))\Big]= 1-\Em_{-x}\Big[\prod_{j=1}^{N_t}(1-g(-x_j(t))\Big],
\enbal
\ee
so that (\ref{20apr1402}) is a special case of (\ref{20apr1406}) with $g(x)=\one(x<0)$. Here, $\Em_x$ refers to a BBM that
starts at~$t=0$ with a single particle at the position $x\in\Rm$ and not at $x=0$. 
The slightly unusual form of~(\ref{20apr1406}) in the right side, with the process
starting at $(-x)$, will be convenient when we look at the limiting measure of the
BBM, also with a flipped sign, as in (\ref{20jun1820}) below.

McKean's interpretation has been used by Bramson~\cite{Bramson1,Bramson2} to establish the following result on the long
term behavior of the solutions to the Fisher-KPP equation. It has been known since the pioneering work of Fisher~\cite{Fisher}
and Kolmogorov, Petrovskii and Piskunov~\cite{KPP} that the Fisher-KPP equation admits traveling wave solutions
of the form  $u(t,x)=U_c(x-ct)$ that satisfy
\begin{equation}\label{20apr1408}
-cU_c'=U_c''+U_c-U_c^2,~~U_c(-\infty)=1,~~U_c(+\infty)=0,
\end{equation}
for all  $c\ge c_*=2$. We will denote by $U(x)$ the traveling wave $U_2(x)$ that corresponds to the minimal speed $c_*=2$:
\begin{equation}\label{jul616}
-2U'=U''+U-U^2,~~U(-\infty)=1,~~U(+\infty)=0.
\end{equation}
Solutions to (\ref{jul616}) are defined up to a translation in space.
We fix the particular translate by requiring that the traveling wave has the asymptotics 
\begin{equation}\label{20mar3114}
U(x)\sim (x+k_0)e^{-x},\hbox{ as $x\to+\infty$},
\end{equation}
with the pre-factor in front of the right side equal to one, and some
$k_0\in\Rm$ that is not, to the best of our knowledge, explicit.  Let now~$u(t,x)$ be the solution
to~(\ref{20apr1402}) with the initial condition~$g(x)$ such that $0\le g(x)\le 1$ for all $x\in\Rm$,
and~$g(x)$ is compactly supported on the right -- there exists~$L_0$ such that $g(x)=0$ for all $x\ge L_0$. It was 
already shown in~\cite{KPP}, for the particular example of~$g(x)=\one(x\le 0)$, 
that there exists a reference frame~$m_{kpp}(t)$ such that
$m_{kpp}(t)/t\to 2$ as $t\to+\infty$ and  
\begin{equation}\label{20apr1412}
u(t,x+m_{kpp}(t))\to U(x)\hbox{ as $t\to+\infty$},
\ee
uniformly on semi-infinite intervals of the form $x\ge K$, for each $K\in\Rm$ fixed. 
Bramson has refined this result, showing that there exists a constant $\hat s[g]$ 
that is known as the Bramson shift corresponding to the initial condition $g$, such that
\begin{equation}\label{20apr1414}
u(t,x+m(t))\to U(x+\hat s[g])\hbox{ as $t\to+\infty$},
\ee
uniformly on semi-infinite intervals of the form $x\ge K$, for each $K\in\Rm$ fixed, with
\begin{equation}\label{20apr1416}
m(t)=2t+\farc{3}{2}\log t.
\ee 
Note that we have chosen the sign of $\hat s[g]$ in (\ref{20apr1414}) in the way that makes the shift
positive for "small" initial conditions that we will consider later. 
In that sense, at a time $t\gg 1$, the solution is  located  at the position $m(t)$. 
A shorter probabilistic proof of this convergence was given recently in~\cite{Roberts}, and PDE proofs 
of various versions on Bramson's results have been obtained in~\cite{HNRR,Lau,NRR1,Uchiyama},
with further refinements in~\cite{Graham,Henderson,NRR2} and especially in the recent
fascinating paper~\cite{BBD}.

\subsection*{The limiting extremal process of BBM and its connection to the Bramson shift }


Motivated by the above discussion, one may 
consider not just the maximal particle but the statistics of the
BBM process re-centered at the location $m(t)$, and ask if it can also
be connected to the solutions to the Fisher-KPP equation.
Let \mbox{$x_1(t)\geq x_2(t)\geq \ldots$} be the positions of the
BBM particles at time $t$, and consider the BBM measure seen from $m(t)$:
\begin{eqnarray}\label{20jun1820}
\cX_t = \sum_{k\leq N_t} \delta_{m(t)-x_k(t)}.
\end{eqnarray}
Recall that $N_t$ is the number of particles alive at time $t$. 
It was shown in~\cite{ABBS,ABK1,ABK2,BD2}   
that there exists a point process $\cX$ so that we have 
\begin{eqnarray}
\cX_t\Rightarrow \cX=\sum_k \delta_{\chi_k}~~\hbox{ as $t\to+\infty$},
\end{eqnarray}
with $\chi_1\leq \chi_2\leq \ldots$, so that $\chi_1$ corresponds to the maximal particle
in the BBM, $\chi_2$ to the second largest, and so on.
In what follows we will call $\cX$ the limiting extremal process
or simply extremal process.  Sometimes it is also called in the literature  the 
decorated Poisson point process, see e.g.~\cite{ABBS}.
The properties of the limit measure $\cX$ are closely related to the 
long time limit of the derivative martingale introduced in~\cite{LS87}:
\begin{equation}\label{20apr1424}
Z_t=\sum_{k\le N_t}(2t-x_k(t))e^{-(2t-x_k(t))}\rightarrow Z~\hbox{ as $t\to+\infty$,~ 
$\P$-a.s.}
\ee
It turns out that there exists a direct connection
	 between $\cX$, $Z$ and the Bramson shift via the Laplace transform of $\cX$.
As explained in Appendix C of~\cite{BD2}, see also~\cite{ABK2}, 
the results of~\cite{LS87} imply that 
for any test function $\psi\geq 0$ that is compactly supported on the right, 
we have the identity
\be \label{20apr2204}
\E\left[ e^{-\cX(\psi)}\right]
=  \E\Big[ e^{-Z e^{-\hat s[\hat\psi]} }\Big], ~~\hat\psi(x)=1-e^{-\psi(x)}.
\ee
Here, for a measure $\mu$ and a function $f$ we use the notation 
\[
\mu(f)=\int_\Rm f(x)\mu(dx).
\]
Note that 
$0\le\hat\psi(x)\le 1$, and $\hat\psi(x)$ is also compactly supported on the right, 
so that its Bramson shift $\hat s[\hat\psi]$ is well-defined and finite. 
In addition to  a special case of (\ref{20apr2204}) 
with~$\hat\psi(x)=\one(x\le 0)$,  
it was shown in~\cite{LS87} that
for each $y\in\Rm$ we have
\begin{eqnarray}
\label{20apr1434}
\E\Big[ e^{-Z e^{-y} }\Big]= 1-U(y).
\end{eqnarray}
This together with (\ref{20apr2204}), characterizes the Laplace transform of $\cX$:
for any test function $\psi\geq 0$ compactly supported on the right, 
we have the duality identity 
\begin{eqnarray}
\label{eq:6_2_1}
 \E\left[ e^{-\cX(\psi)}\right]
= 1- U(\hat s[{\hat\psi}]),~~\hat\psi=1-e^{-\psi}.
\end{eqnarray}
Let us note that the normalization of the traveling wave in (\ref{20mar3114}) implies,
in particular, that  extra shifts appear neither in (\ref{20apr1434}),
nor in (\ref{eq:6_2_1}).  A helpful discussion of the normalization constants
in these identities can be found in Chapter 2 of~\cite{Bovier}. We also explain this
in Section~\ref{sec:probab}.   

It is important to note that, in fact, the  results 
in~\cite{ABBS}, \cite{ABK2} and in   Appendix C of~\cite{BD2} imply the conditional version of~\eqref{20apr2204}, see Lemma~\ref{lem:20_06_1} below:
\be \label{20apr2204a}
\E\left[ e^{-\cX(\psi)}|Z\right]
= e^{-Z e^{-\hat s[\hat\psi]} }.
\ee
In principle, (\ref{20apr2204a}) completely characterizes the
conditional distribution of the  measure $\cX$ in terms of its conditional Laplace
transform. However, the Bramson shift is a very implicit function of the initial condition,
and making the direct use of
(\ref{20apr2204a}) is by no means straightforward. 
{One of the goals of the present paper is precisely to make use of this 
connection to obtain new results about the extremal process $\cX$.}

Let us first illustrate
what kind of results on the Bramson shift we may need on the example of the 
asymptotic growth of $\cX$,
Theorem~1.1 in~\cite{CHL18}, originally
conjectured in~\cite{BD2}. This result says that there exists~$A_0>0$ so that 
\be\label{20apr1426}
\farc{1}{A_0Zxe^x}\cX((-\infty,x])\to1,~~\hbox{ as $x\to+\infty$, in probability.} 
\ee
The proof of (\ref{20apr1426}) 
in~\cite{CHL18} uses purely probabilistic tools. In order to relate this result to 
the Bramson shift and the realm of PDE, 
we can do the following.   
Consider the shifted and rescaled version of the measure~$\cX$, as in (\ref{20apr1426}):
\begin{eqnarray}\label{20apr1432}
Y_n(dx)=  n^{-1} e^{-n} \cX_n(dx),
\end{eqnarray}
where
\begin{eqnarray}\label{20apr1420}
\cX_n=\sum_k \delta_{\chi_k-n},
\end{eqnarray} 
so that
\be\label{20apr1428}
\farc{1}{Zne^n}\cX((-\infty,n])=\farc{1}{Zne^n}\cX_n((-\infty;0])=\farc{1}{Z}Y_n((-\infty;0]).
\ee
We may analyze the conditional on $Z$ Laplace transform of $Y_n$ using
(\ref{20apr2204a}):
given a  non-negative function  $\phi_0(x)$ 
compactly supported on the right, we have 
\begin{equation}
\label{20apr1430}
 \E\left[e^{-Y_n(\phi_0)}|Z\right]=\E\left[ e^{-  n^{-1} e^{-n} \cX_n(\phi_0)}|Z\right]
 =  \E\left[ e^{- \cX(\phi_n)}|Z\right]
 =  \exp\big\{-Z e^{-\hat s[\psi_n]} \big\},
\end{equation} 
with
\begin{equation}\label{20apr1436}
\phi_n(x)= n^{-1} e^{-n}\phi_0(x-n),~~~\psi_n(x)=1-\exp\{-\phi_n(x)\}.
\ee
Note that $\psi_n$ is also compactly supported on the right, so that its Bramson shift is well-defined. 
Furthermore, for $n\gg 1$ the function $\psi_n(x)$ is small: it is of the size $O(n^{-1}e^{-n})$,
as is $\phi_n(x)$.  
Thus,~(\ref{20apr1430}) relates 
the understanding of the conditional on $Z$ weak 
limit of $Y_n$  to the asymptotics of the Bramson 
shift for small initial conditions for the Fisher-KPP equation (\ref{20apr1402}), and this is the
strategy we will exploit in this paper to obtain limit theorems for the process $\cX$.  Let us stress that
a connection between the limiting statistics 
of BBM and the Bramson shift for small initial conditions was already made
in~\cite{BD2}, though with a slightly different objective in mind, and in a different way.

\subsection*{The Bramson shift for small initial conditions}\label{sec:real}

We now state the results for the Bramson shift of the solutions to the Fisher-KPP equation
\begin{equation}\label{jul610}
\pdr{u_\eps}t=\frac{\partial^2 u_{\eps}}{\partial x^2}+u_\eps-u_\eps^2,
\end{equation}
with a small initial condition
\begin{equation}\label{jul612}
u_\eps(0,x)=\eps \phi_0(x),
\end{equation}
that we will need for 
studying the limiting behavior of 
$\cX$. 
Here,~$\eps\ll 1$ is a small parameter, and the function $\phi_0(x)$ is non-negative,
bounded and compactly supported on the right: there exists
$L_0\in\Rm$ such that $\phi_0(x)=0$ for $x\ge L_0$. We will use the notation~$x_\eps=\hat s[\eps\phi_0]$ for the
Bramson shift of $\eps\phi_0$: 
\begin{equation}\label{jul614}
|u_\eps(t,x+m(t)|\to U(x+x_\eps)\to 0\hbox{ as $t\to+\infty$.
uniformly on compact intervals in $x$,}
\end{equation}
We chose the sign of $x_\eps$ in (\ref{jul614}) so that $x_\eps>0$ for 
$\eps>0$ sufficiently small. 
The following proposition gives the asymptotic behavior for $x_\eps$ for small~$\eps>0$
that is sufficiently precise to recover 
(\ref{20apr1426}).
\begin{prop}\label{thm-jun201}
Under the above assumptions on $\phi_0$, we have 
\begin{equation}\label{jun2702}
|x_\eps-\log\eps^{-1}+\log\log\eps^{-1}+\log \bar c|\to 0\hbox{ as $\eps\downarrow 0$,}
\end{equation}
with
\begin{equation}\label{jul620}
\bar c=\farc{1}{\sqrt{4\pi}}\int_{-\infty}^{\infty}e^x\phi_0(x)dx.
\end{equation}
\end{prop}
To formulate the convergence result, 
let $\cC_{c}$ (resp. $\cC^+_{c}$) 
be the space of continuous (resp. non-negative continuous) compactly supported functions 
on $\Rm$, and $\cC_{bc}$ (resp. $\cC^+_{bc}$) 
be the space of bounded continuous (resp. non-negative bounded continuous) 
functions on $\Rm$ compactly  supported on  the right.
Let $\cM_v$ (resp. $\cM_v^+$) be the space of signed (resp. non-negative)
Radon measures on~$\Rm$ such that $|\mu|((-\infty,0])<\infty$, 
equipped with the 
topology $\tau_v$  generated by 
\[
\mu_n \underset{n\rightarrow\infty}{
	\overset{\tau_v}{\longrightarrow}}\mu \iff \mu_n(f)\underset{n\rightarrow\infty}{\longrightarrow} \mu(f) \quad  \forall f\in \calC_{bc}. 
\]
An immediate corollary of Proposition~\ref{thm-jun201} is the following version of
(\ref{20apr1426}). We set
\begin{equation}\label{20mar3102}
\mu(dx)= \frac{1}{\sqrt{4\pi}} e^x \,dx,
\ee
so that 
\be\label{20mar3118}
\bar c= \mu(\phi_0).
\ee
\begin{thm}
\label{lem:28_01_20_1bis} 
Conditionally on $Z$,  we have
\be\label{20mar3022}
Y_n(dx) \underset{n\rightarrow\infty}{\longrightarrow} Z \mu(dx)~~\hbox{{in ${\mathcal M}^+_v$} in probability.}  
\ee
\end{thm}
In other words, $Y_n(dx)$ looks like an exponential shifted by $\log Z$ to the left. 
This  
 theorem also explicitly identifies the constant $A_0$ 
in (\ref{20apr1426}). 
Accordingly, we can reformulate Theorem~\ref{lem:28_01_20_1bis} as follows: 
consider the measures   $ \cX^*_n$ shifted by $\log Z$: 
%
\be\label{20mar3026}
\cX^*_n \equiv \sum_k \delta_{\chi_k+\log Z},
\ee
and
\be\label{20mar3024}
Y^*_n(dx)=  n^{-1} e^{-n} \cX^*_n(dx). 
\ee
This gives the following version of Theorem~\ref{lem:28_01_20_1bis}:
\begin{cor} 
\label{lem:28_01_20_2bis}
We have 
\be\label{20mar3028}
Y^*_n(dx)\underset{n\rightarrow\infty}{\longrightarrow}   \mu(dx), 
~~\hbox{{in ${\mathcal M}^+_v$} in probability.}  
\ee
\end{cor}

As we have mentioned, Theorem~\ref{lem:28_01_20_1bis} and Corollary~\ref{lem:28_01_20_2bis}
are not really new, except for identifying the constant $A_0$, even though the approach via the asymtptotics
of the Bramson shift produces
an analytic rather than a probabilistic proof.
In order to obtain 
genuinely new results on the fluctuations of $Y_n(dx)$ around $Z\mu(dx)$, we will need a finer 
asymptotics for the shift $x_\eps$ than in 
Proposition~\ref{thm-jun201}. 
Let us define the constant
\begin{equation}\label{20apr104}
\begin{aligned}
\bar c_1&=\farc{1}{\sqrt{4\pi}}\int_{-\infty}^\infty x e^x \phi_0(x)dx,
\enbal
\ee
that depends on the initial condition $\phi_0$, as does $\bar c$ in (\ref{jul620}), and 
universal constants 
\be\label{20apr1440} 
g_\infty= \int_0^1 e^{z^2/4}\int_z^\infty e^{-y^2/4}dy dz
-2\int_1^\infty e^{z^2/4}\int_z^\infty \farc{1}{y^2}e^{-y^2/4}dy,
\end{equation}
and
\begin{equation}\label{20apr102}
m_1=\farc{3}{2}g_\infty+k_0+\farc{1}{2},
\end{equation}
that do not depend on $\phi_0$. Here, $k_0$ is the constant that appears in 
the asymptotics (\ref{20mar3114}) for $U(x)$. 
The following theorem allows us to  
obtain convergence in law of the fluctuations of $Y_n$. 
\begin{thm}\label{conj-jun1102-19}
Under the above assumptions on $\phi_0$, we have the asymptotics
\begin{equation}\label{19jun1102}
x_\eps=\log\eps^{-1}-\log\log\eps^{-1}-\log \bar c
-\frac{2\log\log\eps^{-1}}{\log\eps^{-1}}-\Big(m_1-\log\bar c+\farc{\bar c_1}{\bar c}\Big)
\frac{1}{\log\eps^{-1}}+
O\Big(\frac{1}{(\log\eps^{-1})^{1+\gamma}}\Big),
\end{equation}
as $\eps\downarrow 0$, with some $\gamma>0$. 
\end{thm}

The first two terms in (\ref{jun2702}) and
(\ref{19jun1102}) have been predicted in \cite{BD2} in addressing a different BBM question,
using an informal Tauberian type argument that we were not able to make rigorous. 
The rest of terms have not been predicted, to the best of our knowledge.
The proof in the current paper 
does not seem to be directly related  to the arguments of~\cite{BD2} 
but the general approach to the statistics of BBM via the Bramson shift asymptotics for small initial conditions
comes from~\cite{BD2}. 

\subsection*{{Weak convergence of the fluctuations of the extremal  process}}
 
{Theorem~\ref{lem:28_01_20_1bis} indicates that to
obtain the limiting behavior of the fluctuations of $Y_n(dx)$} 
one should consider a rescaling of the signed measure  
\be\label{20apr1602}
Y_n(dx)-Z\mu(dx)=Y_n(dx)-\farc{1}{\sqrt{4\pi}}e^{(x+\log Z)}dx.
\ee
It turns out, however, that there is an extra small deterministic shift of the
exponential profile that needs to be performed before the rescaling: a better
object to rescale is not as in (\ref{20apr1602}) but 
\be\label{20apr1604}
Y_n(dx)-\farc{1}{\sqrt{4\pi}}e^{(x+\log Z+e_n)}dx,
\ee
with an extra deterministic correction 
\be\label{20apr1606}
e_n=\log\Big(1+\farc{2\log n}{n}\Big)\approx \frac{2\log n}{n}.
\ee
The properly rescaled measures are
\begin{equation}\label{20apr1610}
\bal
&V_n(\phi_0) = n \Big(Y_n(\phi_0) - (1+\frac{2\log n}{n})Z\mu(\phi_0)\Big), 
\\
&V^*_n(\phi_0) = n \Big(Y^*_n(\phi_0) - (1+\frac{2\log n}{n})\mu(\phi_0)\Big).
\enbal
\end{equation}
Let us set
\be\label{20apr1608}
\nu(dx) = (4\pi)^{-1/2} xe^x \,dx,
\ee
%
and let $\{ R_t, t\geq 0\}$ be
a spectrally positive $1$-stable stochastic process  with the Laplace transform 
$$
\E\left[e^{- \lambda R_t}\right]=e^{t\lambda\log\lambda},\quad \forall \lambda >0, \;t\geq 0,
$$
such that $\{R_t, t\geq 0\}$ is independent of $Z$.
The main probabilistic result of this paper is the following   
theorem describing the limiting behavior of fluctuations of the extremal process.
We use the notation~$\Rightarrow$ for the convergence in distribution.
\begin{thm}
\label{prop:1}
(i) 
Conditionally on $Z$, we have 
\begin{eqnarray}
\label{eq:27_06_5}
V_n \Rightarrow \bbL_{Z}(dx)\quad\hbox{in ${\mathcal M}_v$, as $n\rightarrow\infty$.}
\end{eqnarray}
Here, $\bbL_{Z}$ is a random measure such that 
\begin{equation}\label{20mar3112}
\bbL_{Z}(dx) = R_{Z} \mu(dx) + Z (m_1 \mu(dx)+\nu(dx)),
\end{equation}
and $m_1$ is the constant that appears in (\ref{19jun1102}).\\
(ii) 
We also have
\begin{eqnarray}
\label{eq:27_06_5b}
V^*_n(dx)  \Rightarrow \bbL_{1}(dx) \quad\hbox{in ${\mathcal M}_v$, as $n\rightarrow\infty$.}
\end{eqnarray} 
Here, $\bbL_1(dx)$ is   a random measure such that
\be\label{20mar3114bis}
\bbL_1(dx) = R_{1} \mu(dx) + (m_1 \mu(dx)+\nu(dx)).  
\ee
\end{thm}
One can immediately deduce from the above theorem, that, conditionally on $Z$, 
for  any test function~$\phi_0\in \calC^+_{bc}$,~$\bbL_Z(\phi_0)$ is a spectrally positive $1$-stable random variable with the Laplace transform:
\begin{eqnarray}\label{20mar3138}
\E \left[e^{- \lambda \bbL_{Z}(\phi_0) }\big|Z\right] =
e^{Z\mu(\phi_0)\lambda\log \lambda +\lambda Z\mu(\phi_0)\log \mu(\phi_0) 
-\lambda{Z}(m_1\mu(\phi_0)+\nu(\phi_0))  },\;\forall \lambda>0,
\end{eqnarray}


We note that, in the study of BBM, both the $1$-stable process behavior 
and a small deterministic correction
have appeared in a related but different context 
for the convergence~(\ref{20apr1424}) 
of $Z_t$ to $Z$ as~$t\to+\infty$ in~\cite{MP}.  There, the correction is of the 
form $\log t/\sqrt{t}$, and is close in spirit to $e_n\sim\log n/n$ due to the diffusive
nature of the
space-time scaling, though the exact translation of the corrections appearing
here and in~\cite{MP} is
not quite clear. The $\log t/t$ correction to the front location has also been
observed in~\cite{BBD} and proved in~\cite{Graham}. 
The $1$-stable like tails for $Z$ itself have been observed already
in~\cite{BerBerSch} for the BBM, and in~\cite{Madaule} for the branching random walk. 

The paper is organized is follows. Section~\ref{sec:probab} explains how the 
the probabilistic statements, Theorem~\ref{lem:28_01_20_1bis}
and Theorem~\ref{prop:1}, follow from the correspomnding asymptoptics for the Bramson
shift~$x_\eps$ in Proposition~\ref{thm-jun201} and Theorem~\ref{conj-jun1102-19},
as well as the  the duality identities (\ref{20apr1434}) and
(\ref{eq:6_2_1}). The rest of the paper is devoted to the proof of
Proposition~\ref{thm-jun201} and Theorem~\ref{conj-jun1102-19}. 
Section~\ref{sec:self-similar} introduces the self-similar variables that 
are used throughout the proofs, and reformulates the required results as
Propositions~\ref{prop-jun2504} and~\ref{conj-jun1104-19}. We also explain
in that section how the values of the specific constants that appear in the asymptotics
for $x_\eps$ come about. 
Section~\ref{sec:linear} contains the analysis of the linear Dirichlet problem
that appears throughout the PDE approach to the Bramson shift~\cite{HNRR,NRR1,NRR2}.
More specifically, we describe an approximate solution to the adjoint linear problem that is 
one reason for the logarithmic correction in $e_n$ in (\ref{20apr1606}).
We should stress  that this is not the only contribution to 
the $\log n/n$ term -- the second one comes from the nonlinear term in the
Fisher-KPP equation. Section~\ref{sec:prop-lln} contains the proof of 
Proposition~\ref{prop-jun2504}. Of course, this Proposition  is just a weaker version
of Proposition~\ref{conj-jun1104-19}, in the same vein as 
Proposition~\ref{thm-jun201} is an immediate consequence of Theorem~\ref{conj-jun1102-19}.
However, the proof of Proposition~\ref{prop-jun2504} is much shorter than that of 
Proposition~\ref{conj-jun1104-19}, so we present it separately for the convenience of the reader.
Section~\ref{sec:5} contains the proof of Proposition~\ref{conj-jun1104-19},
with some intermediate steps proved in Section~\ref{sec:aux}. Finally,
Appendix~\ref{sec:appendix} contains an auxiliary lemma on the continuity of
the Bramson shift.

We use the notation $C$, $C'$, $C_0$, $C_0'$, etc. for various constants that do not depend on $\eps$,
and can change from line to line.

{\bf Acknowledgment.}  We are deeply indebted to Lisa Hartung for illuminating discussions
during the early stage of this work, and 
are extremely grateful to \'Eric Brunet and Julien Berestycki for generously sharing 
their deep understanding of various aspects
of the BBM and the Bramson shift, without which this work would not be possible. 
The work of JMR was  supported
by from the European Research Council under the European Union’s Seventh Framework Programme (FP/2007-2013) ERC Grant Agreement n. 321186 - ReaDi, 
and from the ANR NONLOCAL project (ANR-14-CE25-0013). The work or LM and LR was supported by
a US-Israel BSF grant, 
LR was partially supported by NSF grants DMS-1613603 and DMS-1910023, and ONR grant N00014-17-1-2145. 

\section{The proof of the probabilistic statements}\label{sec:probab}

In this section, we explain how the probabilistic results, 
Theorem~\ref{lem:28_01_20_1bis} and Theorem~\ref{prop:1},
follow from Proposition~\ref{thm-jun201} and Theorem~\ref{conj-jun1102-19}, respectively.
 
\subsection*{The duality identity}

We first briefly explain the duality identities (\ref{20apr1434}) and
(\ref{eq:6_2_1}), and, in particular, the fact that
no extra shift of the wave is needed in these relations, 
as soon as the wave normalization (\ref{20mar3114})
is fixed. To see this, note that in the formal
limit~$\psi(x)\to \infty\cdot\one(x\le y)$, 
we have $\hat\psi(x)=\one(x\le y)$, and~(\ref{eq:6_2_1})
reduces it to
\be\label{20apr2302}
\Pm[\chi_1\le y]=1-U(\hat s[\chi_0]+y),~~\chi_0(x)=\one(x\le 0). 
\ee
This is simply the definition of the Bramson shift, combined with the probabilistic
interpretation~(\ref{20apr1404}) of the solution to (\ref{20apr1402}) with the initial 
condition $u(0,x)=\one(x\le 0)$. Thus, no extra shift is needed neither in
(\ref{20apr2302}), nor in~(\ref{eq:6_2_1}). 
The generalization of (\ref{20apr2302}) to (\ref{eq:6_2_1}) is explained in Appendix C
of~\cite{BD2}. Some normalization constants appear in the discussion there, but as we see,
they are not needed  once the wave is normalized by (\ref{20mar3114}) rather than by 
\begin{equation}\label{20apr2308}
\int_{-\infty}^\infty x\tilde U'(x)dx=0,
\end{equation}
as in \cite{BD2}. 

As for the derivative martingale identity (\ref{20apr1434}), it is simply an immediate 
consequence of expressions~(10) and (11) in~\cite{LS87}, 
combined with  the normalization (\ref{20mar3114})
that implies that $C=1$ in both of these equations in~\cite{LS87}.

\subsection*{{Proof of Theorem~\ref{lem:28_01_20_1bis}}}

We now prove Theorem~\ref{lem:28_01_20_1bis}. 
For an arbitrary $\phi_0\in\cC^+_{bc}$, we set, as in (\ref{20apr1436}),
\begin{equation}\label{20mar3016}
\phi_n(x)= n^{-1} e^{-n}\phi_0(x-n),~~~ \tilde\phi_n(x)= n^{-1} e^{-n}\phi_0(x),
\end{equation}
and 
\be\label{20mar3018}
\psi_n(x)=1-\exp\{-\phi_n(x)\}, ~~\tilde\psi_n(x)=1-\exp\{-\tilde\phi_n(x)\}, 
\ee
with 
\be\label{20mar3020}
\psi_0(x)=1-\exp\{-\phi_0(x)\}.
\ee
Note that  both $\psi_n$ and $\tilde\psi_n$ look like "small step" initial conditions:
\[
\psi_n(x)\approx \phi_n(x),~~\tilde\psi_n(x)\approx  \tilde\phi_n(x)\hbox{ as $n\to+\infty$}. 
\]  
We need the following result  that, in fact, follows easily from \cite{ABBS}, \cite{ABK2} and Appendix~C in~\cite{BD2}. 
We provide the proof for the sake of completeness we provide 
the proof. 
\begin{lemma}\label{lem:20_06_1}
For any $\phi_0\in \cC^+_{bc}$, we have
\be
\label{eq:27_01_1}
\E\Big[ e^{-\cX(\phi_0)} |Z \Big]
= \exp\big\{-Z e^{-\hat s[\psi_0]}\big\},
\ee
and hence 
\begin{eqnarray}
\label{eq:20_6_1}
\E\Big[ e^{-Y_n(\phi_0)} |Z \Big]
= \exp\big\{-Z e^{-\hat s[\psi_n]}\big\}.
\end{eqnarray} 
\end{lemma}
{\bf Proof.}
Fix an arbitrary $\phi_0\in \cC^+_{c}$ and let $\{\cF_t\}_{t\geq 0}$ 
be a  filtration generated by $\{\cX_t\}_{t\geq 0}$.
Then, by the definition (\ref{20jun1820}) of $\{\cX_t\}_{t\geq 0}$, 
the Markov property and \eqref{20apr1406}, we get (with some ambiguity of notation we set $\E=\E_0$: the expectation for BBMs starting with one particle at $0$)
\be 
\bal
\E\Big[ e^{-\cX_t(\phi_0)}&| \cF_s\Big]=
\E\Big[\prod_{j=1}^{N_t}e^{-\phi_0(2t-(3/2)\log t-x_j(t)}| \cF_s
\Big]\\
&=
\prod_{j=1}^{N_s} \Em_{x_j(s)}\Big[\prod_{i=1}^{N_{t-s}}
\Big(1-(1- e^{-\phi_0(2(t-s)-({3}/{2})\log(t-s)-x_i(t-s)+2s+
({3}/{2})\log((t-{s})/{t}))}\Big)\Big]
\\
\label{eq:28_05_1}
&= 
\prod_{j=1}^{N_s} (1-u(t-s, m(t-s)-x_j(s)+2s+\frac{3}{2}\log(1-\frac{s}{t}))),\quad  \P-\text{a.s.}
\enbal
\ee
Here, $u(t,x)$ is the solution  
to \eqref{20apr1402} with the initial condition $u(0,x)=\psi_0(x)=1-\exp(-\phi_0(x))$. 
Next, we use \eqref{20apr1414} to get 
\be 
\bal
\lim_{t\rightarrow\infty}
&\prod_{j=1}^{N_s} (1-u(t-s, m(t-s)-x_j(s)+2s+\frac{3}{2}\log(1-\frac{s}{t})))
= 
\prod_{j=1}^{N_s} (1-U(\hat s[\psi_0]-x_j(s)+2s))\\
\label{eq:28_05_2}
&=
\exp\Big\{-\sum_{j=1}^{N_s} - \log(1-U(\hat s[\psi_0]-x_j(s)+2s)\Big\}.
\enbal
\ee
Following the derivations in (2.4.9)-(2.4.12)  
in~\cite{Bovier} (see also (23)-(25) in~\cite{LS87}) and using~\eqref{20mar3114} which 
gives $C=1$ in the above references we obtain
\be 
 \lim_{s\rightarrow\infty} \exp\Big\{-\sum_{j=1}^{N_s} - \log(1-U(\hat s[\psi_0]-x_j(s)+2s)\Big\}
= 
 \label{eq:28_05_3}
\exp\Big\{-Z e^{-\hat s[\psi_0]}\Big\},\quad \P-\text~{\rm a.s.}
\ee
Fix an arbitrary bounded continuous function $h$.  
Putting together (\ref{eq:28_05_1}), (\ref{eq:28_05_2}), and  (\ref{eq:28_05_3}), we get,
\be 
\bal
\lim_{s\rightarrow\infty}  &\lim_{t\rightarrow\infty} 	\E\Big[ h(Z_s)e^{-\cX_t(\phi_0)}\Big]
 	= \lim_{s\rightarrow\infty} \lim_{t\rightarrow\infty}
 	\E\Big[ h(Z_s) \E\Big[e^{-\cX_t(\phi_0)}| \cF_s\Big]\Big]
 	\\
 	&= \lim_{s\rightarrow\infty}
 	\E\Big[ h(Z_s) e^{-\sum_{j=1}^{N_s} - \log(1-U(\hat s[\psi_0]-x_j(s)+2s)}\Big]
 	\label{eq:17_06_1}
 	=
 	\E\Big[ h(Z)e^{-Z e^{-\hat s[\psi_0]}}\Big].
\enbal
\ee
We used 
the bounded convergence theorem in the last equality.
 
Recall, that by Theorem~2.1 in~\cite{ABBS}, the pair $(\bar\cX_t, Z_t)$ converges 
in distribution to $(\bar\cX, Z)$ as~$t\rightarrow\infty$, where $\bar\cX_t$ 
 (resp. $\bar\cX$) is just the measure $\cX_t$ (resp. $\cX$) shifted by $\log Z$.
This together with a.s. convergence of $Z_t$ to $Z$ implies also convergence in distribution of $(\cX_t, Z_t)$  to $(\cX, Z)$ as $t\rightarrow\infty$,
and thus
\begin{eqnarray*}
\E\Big[ h(Z)e^{-\cX(\phi_0)}\Big]
 = 	 \lim_{t\rightarrow\infty} 	\E\Big[ h(Z_t)e^{-\cX_t(\phi_0)}\Big].
 \end{eqnarray*}
From this we get
 \[
 \bal
 &\left| \E\Big[ h(Z)e^{-\cX(\phi_0)}\Big]-\E\Big[ h(Z)e^{-Z e^{-\hat s[\psi_0]}}\Big]\right|
 = 	 	\left| \lim_{t\rightarrow\infty} 	\E\Big[ h(Z_t)e^{-\cX_t(\phi_0)}\Big]-\E\Big[ h(Z)e^{-Z e^{-\hat s[\psi_0]}}\Big]\right|\\
 &\leq 
 \left| \lim_{t\rightarrow\infty} 	\E\Big[ (h(Z_t)-h(Z_s))e^{-\cX_t(\phi_0)}\Big]\right|
 +  \left| \lim_{t\rightarrow\infty} 	\E\Big[ h(Z_s)e^{-\cX_t(\phi_0)}\Big]-\E\Big[ h(Z)e^{-Z e^{-\hat s[\psi_0]}}\Big]\right|,  \; \forall s>0. 
 \enbal
 \]
 Therefore, we have 
\[
\bal
\left| \E\Big[ h(Z)e^{-\cX(\phi_0)}\Big]-\E\Big[ h(Z)e^{-Z e^{-\hat s[\psi_0]}}\Big]\right|
 	&\leq \limsup_{s\rightarrow \infty}
  \lim_{t\rightarrow\infty} 	\E\Big[ 	\left|h(Z_t)-h(Z_s)\right|\Big]\\
 	&+  \left| \limsup_{s\rightarrow \infty}\lim_{t\rightarrow\infty} 	\E\Big[ h(Z_s)e^{-\cX_t(\phi_0)}\Big]-\E\Big[ h(Z)e^{-Z e^{-\hat s[\psi_0]}}\Big]\right| =0,
\enbal
\]
where convergence to zero of the first term follows from a.s. convergence of $Z_t$ to $Z$ and 
the bounded convergence theorem. The second term converges to zero by (\ref{eq:17_06_1}).
Hence, we have
\[
	\E\Big[ h(Z)e^{-\cX(\phi_0)}\Big]=\E\Big[ h(Z)e^{-Z e^{-\hat s[\psi_0]}}\Big],
\]
for any $\phi_0\in \cC^+_{c}$ and any bounded continuous function $h$, which implies (\ref{eq:27_01_1}) for any $\phi_0\in \cC^+_c$. The
extension to
$\phi_0\in \cC^+_{bc}$ follows  via approximation and  a kind of continuity of $\hat s[\psi]$ 
in functions in $\cC^+_{bc}$ and bounded by $1$,
made precise in Lemma~\ref{lem-20june1802} in Appendix~\ref{sec:appendix}. \qed

We continue the proof of Theorem~\ref{lem:28_01_20_1bis}. Recall that we fixed an arbitrary $\phi_0\in \cC^+_{bc}$.
It follows from Proposition~\ref{thm-jun201}, with $\eps_n=n^{-1}\exp(-n)$ that 
the Bramson shift appearing in the right side of~(\ref{eq:20_6_1}) has the asymptotics
\begin{equation}\label{eq:28_01_20_2}
\hat s[\psi_n]=-n+ \hat s[\tilde\psi_n]=-n +\log(ne^n)-\log(n+\log n)-\log\bar c+o(1)
\to -\log\mu(\phi_0), \hbox{ as $n\to+\infty$,}
\end{equation} 
with the measure $\mu$ defined in (\ref{20mar3102}), and $\tilde\psi_n$ as in (\ref{20mar3018}). 

To  get the weak convergence of the measures $Y_n(dx)$, it is sufficient to get the convergence of~$Y_n(\phi_0)$ for any $\phi_0\in\cC^+_{bc}$. 
Thus, we will check convergence of the Laplace transforms of $Y_n(\phi_0)$ (conditionally on~$Z$).
As a consequence of 
Lemma~\ref{lem:20_06_1}  and (\ref{eq:28_01_20_2}), we obtain
\begin{equation}
\label{eq:15_5_8}
\lim_{n\rightarrow \infty} \E\left[ e^{-  Y_n(\phi_0)}|Z\right]
=\lim_{n\rightarrow \infty} e^{-Z e^{-\hat s[\psi_n]}}=e^{-Z\mu(\phi_0)}.
\end{equation} 
Since $\phi_0\in\cC_{bc}^+$  was arbitrary, this  implies  
that, conditionally on $Z$, 
\be
Y_n\Rightarrow Z\mu,\quad \hbox{in ${\mathcal M}^+_v$} ,\;\mbox{as}\; n\rightarrow \infty.
\ee
Therefore, conditionally on $Z$,  we have
\be\label{20mar3022bis}
Y_n(dx)\rightarrow Z \mu(dx)~~\hbox{in ${\mathcal M}^+_v$} ,\;\mbox{as}\; n\rightarrow \infty, \hbox{in probability,}  
\ee
finishing the proof of Theorem~\ref{lem:28_01_20_1bis}.~$\Box$

Corollary~\ref{lem:28_01_20_2bis} is an immediate consequence 
of Theorem~\ref{lem:28_01_20_1bis}.

\subsection*{ Proof of Theorem~\ref{prop:1}}

We will prove only part (i)   since the proof of  (ii) goes along the same lines.
Once again, to get the weak convergence of measures $V_n(dx)$, it is sufficient to get the convergence of $V_n(\phi_0)$ for any~$\phi_0\in\cC^+_{bc}$, and 
we will check the convergence of the Laplace transforms of $V_n(\phi_0)$ (conditionally on $Z$).
Fix an arbitrary $\phi_0\in\cC^+_{bc}$. We have
\begin{equation}\label{20mar3130}
\begin{aligned}
\E\Big[e^{-V_n(\phi_0)}\big| Z\Big] &= 
\E\Big[\exp\Big(-n \Big(Y_n(\phi_0) - (1+\frac{2\log n}{n})Z\mu(\phi_0)\Big)\Big)\big| Z\Big]\\
&=\E\Big[\exp\big(-nY_n(\phi_0)\big)\big|Z\Big]\exp\Big[\big(n+2\log n)Z\mu(\phi_0)\Big]\\
&= \exp\Big[-Z e^{-\hat s[\eta_n]}\Big]\exp\Big[\big(n+2\log n)Z\bar c\Big].
\end{aligned}
\end{equation}
We used expression (\ref{20mar3118}) for $\bar c$ and identity (\ref{eq:20_6_1}) above, with $\phi_0$ replaced by $n\phi_0$, and
\[
\eta_n(x)=1-e^{-e^{-n}\phi_0(x-n)}.
\]
Let us also introduce  
\[
\zeta_n(x)=n\phi_n(x)=e^{-n}\phi_0(x-n),~~
x_n= \hat s[{\zeta_n}],~~ 
\tilde x_n = \hat s[\eta_n],
\]
so that (\ref{20mar3130}) can be written as
\begin{equation}\label{20mar3132}
\begin{aligned}
\E\Big[e^{-V_n(\phi_0)}\big| Z\Big] &= \exp\Big[-Z\big(e^{-\tilde x_n} -\big(n+2\log n)\bar c\big)\Big].
\end{aligned}
\end{equation}
Note that by Theorem~\ref{conj-jun1102-19} with $\eps_n=e^{-n}$ we get 
\begin{equation}\label{20mar3120}
\begin{aligned}
\tilde x_n &=  x_n + O(e^{-n} \left\| \phi_0\right\|_{\infty})=-n
+\log\eps_n^{-1}-\log\log\eps_n^{-1}-\log \bar c
-\frac{2\log\log\eps_n^{-1}}{\log\eps_n^{-1}}\\
&-\Big(m_1-\log\bar c+\farc{\bar c_1}{\bar c}\Big)
\frac{1}{\log\eps_n^{-1}}+
O\Big(\frac{1}{(\log\eps_n^{-1})^{1+\delta}}\Big) \\
&=-\log n - \log\bar c -2\frac{\log n}{n} 
-\farc{1}{n}\Big(m_1-\log\bar c+\frac{\bar c_1}{\bar c}\Big)+O(n^{-1-\delta}).
\end{aligned}
\end{equation}
Using this in (\ref{20mar3132}) gives
\begin{equation}\label{20mar3133}
\begin{aligned}
\E\Big[e^{-V_n(\phi_0)}\big| Z\Big] &= \exp\Big[-Z\big(e^{-\tilde x_n} -\bar cn-2\bar c\log n\big)\Big]\\
&= \exp\Big[-\bar c Zn\Big\{\exp\Big[2\frac{\log n}{n}+\farc{1}{n}\Big(m_1-\log\bar c+\frac{\bar c_1}{\bar c}\Big)+O(n^{-1-\delta})\Big]
-1-\farc{2\log n}{n}\Big\}\Big]\\
&=\exp\Big[-\bar c Zn\Big\{\farc{1}{n}\Big(m_1-\log\bar c+\frac{\bar c_1}{\bar c}\Big)+O(n^{-1-\delta})\Big\}\Big]\\
&\to \exp\Big[-\bar c Z \Big(m_1-\log\bar c+\frac{\bar c_1}{\bar c}\Big) \Big]=\exp\Big[Z\Big(\bar c\log\bar c-\bar c m_1-\bar c_1\Big)\Big],
~~\hbox{ as $n\to+\infty$.}
\end{aligned}
\end{equation}
By the definition of $\bbL_1$, $\bar c$ and $\bar c_1$ we get
\begin{eqnarray*}
\E\Big[e^{-\bbL_1(\phi_0)}\big| Z\Big]= \exp\Big[Z\bar c\log\bar c -  Z(\bar c m_1 + \bar c_1)\Big],
\end{eqnarray*}
and since $\phi_0\in\cC^+_{bc}$ was arbitrary we are done. \qed 


\section{The solution asymptotics in self-similar variables} \label{sec:self-similar}

Proposition \ref{thm-jun201} is a consequence of the following two steps. 
The first result connects the Bramson shift of a solution to the Fisher-KPP
equation with a small initial condition to the asymptotics of the solution 
to a problem
in the self-similar variables with an initial condition shifted far to the right. 
\begin{prop}\label{prop-jun2502}
Let $r_\ell$ be the solution to 
\begin{eqnarray}\label{jun20104}
\pdr{r_\ell}{\tau}-\frac{\eta}{2}\pdr{r_\ell}{\eta}
-\pdrr{r_\ell}{\eta}- r_\ell
+\farc{3}{2}e^{-\tau/2}\pdr{r_\ell}{\eta}
+ e^{3\tau/2-\eta\exp(\tau/2)}r_\ell^2=0,~~\tau>0,~~\eta\in\Rm,
\end{eqnarray}
with the initial condition
$r_\ell(0,\eta)=\psi_0(\eta-\ell)$,
where $\psi_0(\eta)=e^\eta\phi_0(\eta)$. 
Then, for each $\ell>0$, the function~$r_\ell(\tau,\eta)$ has the asymptotics
\begin{equation}\label{may1012}
r_\ell(\tau,\eta)\sim
r_\infty(\ell)\eta e^{-\eta^2/4},\hbox{ as $\tau\to+\infty$, for $\eta>0$}.
\end{equation}
Furthermore, the Bramson shift that appears in Proposition~\ref{thm-jun201} is given by
\begin{equation}\label{19jun1210}
x_\eps=\log\eps^{-1}-\log r_\infty(\ell_\eps), \hbox{ with $\ell_\eps=\log\eps^{-1}$.}
\end{equation}
\end{prop}
The second result, at the core of the proof of Proposition~\ref{thm-jun201}, describes the
asymptotics of $r_\infty(\ell)$ for large $\ell$.
\begin{prop}\label{prop-jun2504}
The function $r_\infty(\ell)$ satisfies the following asymptotics: 
\begin{equation}\label{jun2502}
r_\infty(\ell)=\bar c\ell+O(\log\ell),~~\hbox{ as $\ell\to+\infty$},
\end{equation}
with the constant $\bar c$ as in (\ref{jul620}). 
\end{prop} 
To prove Theorem~\ref{conj-jun1102-19}, we  refine
Proposition~\ref{prop-jun2504} to the following. 
\begin{prop}\label{conj-jun1104-19}
The function $r_\infty(\ell)$ satisfies the following asymptotics: 
\begin{equation}\label{jan810bis}
r_\infty(\ell)=\bar c\ell+2\bar c \log\ell
+m_1\bar c +\bar c_1-\bar c\log\bar c+ O(\ell^{-\delta}),
\end{equation}
%
with the constants $\bar c$, $\bar c_1$ and $m_1$ as in (\ref{jul620}),
(\ref{20apr104}) and (\ref{20apr102}).
\end{prop}

Using (\ref{19jun1210}), we obtain from Proposition~\ref{conj-jun1104-19} that
\begin{equation}\label{19jun1208}
\begin{aligned}
x_\eps&=\ell_\eps-\log r_\infty(\ell_\eps)=
\ell_\eps-\log\Big(\bar c\ell_\eps+2\bar c\log\ell_\eps+m_1\bar c-
\bar c\log\bar c+\bar c_1 +O(\ell_\eps^{-\delta})\Big)\\
&=\ell_\eps-\log\ell_\eps-\log \bar c- 
2 \farc{\log\ell_\eps}{\ell_\eps}-\frac{m_1}{\ell_\eps}+\farc{\log\bar c}{\ell_\eps}-
\frac{\bar c_1}{\bar c\ell_\eps}
+O(\ell_\eps^{-1-\delta}),
\end{aligned}
\end{equation}
which proves Theorem~\ref{conj-jun1102-19}.
Thus, our goal is
to prove Proposition~\ref{conj-jun1104-19}. 

Of course, Proposition~\ref{prop-jun2504} in an immediate consequence
of Proposition~\ref{conj-jun1104-19}. However, as its proof is both much shorter
and helpful in the proof of the latter, we present its proof separately 
in Section~\ref{sec:prop-lln}. 
 
\subsection*{Common sense scaling arguments}

In order to verify that the constants in Proposition~\ref{conj-jun1104-19}
are plausible, let us see assume that we have the asymptotics
\begin{equation}\label{19jun1108}
r_\infty(\ell)=\bar c\ell+m_0\bar c\log\ell+m_1\bar c+m_2 \bar c\log\bar c+m_3\bar c_1+
O(\ell^{-\delta}),~~\hbox{ as $\ell\to+\infty$},
\end{equation}
and see what simple arguments say about the possible
values of the coefficients $m_0$, $m_1$, $m_2$ and~$m_3$.
First, consider a shifted initial condition
$\phi_0^L(x)=\phi_0(x-L)$.
Then, the shift $x_\eps$ given 
by (\ref{19jun1208}) 
should also change by $L$, so that
\begin{equation}\label{jan2104}
x_\eps^L=x_\eps -L.
\end{equation}
Note that
\begin{equation*}
\begin{aligned}
&\bar c_L=\farc{1}{\sqrt{4\pi}}\int_{-\infty}^{\infty}e^x\phi_0(x-L)dx=e^L\bar c,\\
&\bar c_1^L=\farc{1}{\sqrt{4\pi}}\int_{-\infty}^\infty x e^x\phi_0(x-L)dx=e^L\bar c_1+
Le^L\bar c.
\end{aligned}
\end{equation*}
Using (\ref{19jun1108})  gives
\begin{equation*}
\begin{aligned}
x_\eps^L&=\log\eps^{-1}-\log r_\infty^L(\ell_\eps)\\
&=\log\eps^{-1}-\log\log\eps^{-1}-\log \bar c_L- 
m_0\farc{\log\log\eps^{-1}}{\log\eps^{-1}}-\frac{m_1}{\log\eps^{-1}}-
\frac{m_2\log\bar c_L}{\log\eps^{-1}} -\farc{m_3\bar c_1^L}{\bar c_L\log\eps^{-1}}
+O(\ell_\eps^{-1-\delta})\\
&=\log\eps^{-1}-\log\log\eps^{-1}-\log \bar c-L- 
m_0\farc{\log\log\eps^{-1}}{\log\eps^{-1}}-\frac{m_1}{\log\eps^{-1}}-
\frac{m_2\log\bar c}{\log\eps^{-1}}-\frac{m_2L }{\log\eps^{-1}}\\
&-m_3\farc{e^L\bar c_1+Le^L\bar c}{\bar c e^L\log\eps^{-1}}+O(\ell_\eps^{-1-\delta})
=x_\eps-L-\frac{m_2L }{\log\eps^{-1}}
- \farc{m_3L}{\log\eps^{-1}}+O(\ell_\eps^{-1-\delta}).
\end{aligned}
\end{equation*}
This means that for (\ref{jan2104}) to hold we must have
\begin{equation}\label{20jan1916}
m_3=-m_2.
\end{equation}
Note that in (\ref{jan810bis}) we have $m_2=-1$ and $m_3=1$, so that (\ref{20jan1916})
holds. 
 
The second invariance is to consider an initial condition   
$\phi_0^\lambda(x)=\lambda\phi_0$. This is equivalent to
keeping $\phi_0$ intact and replacing $\eps$ by $\eps_\lambda=\eps\lambda$. 
If we replace $\phi_0$ by $\lambda\phi_0$ in (\ref{19jun1208}) and keep $\eps$ unchanged, this corresponds to
replacing $\bar c$ by $\lambda \bar c$ and $\bar c_1$ by $\lambda\bar c_1$, which gives
\begin{equation}\label{20jan1918}
\begin{aligned}
x_\eps^\lambda&=\log\eps^{-1}-\log\log\eps^{-1}-\log \bar c- \log\lambda-
m_0\farc{\log\log\eps^{-1}}{\log\eps^{-1}}-\frac{m_1}{\log\eps^{-1}}-
\frac{m_2\log\bar c}{\log\eps^{-1}}-\frac{m_2\log\lambda}{\log\eps^{-1}}\\
&-m_3\farc{\bar c_1}{\bar c\log\eps^{-1}}
+O(\ell_\eps^{-1-\delta})=x_\eps- \log\lambda-\frac{m_2\log\lambda}{\log\eps^{-1}}
+O(\ell_\eps^{-1-\delta}).
\end{aligned}
\end{equation}
If, instead, we replace $\eps$ by $\eps_\lambda=\eps\lambda$ and keep $\phi_0$ intact, 
we get  from (\ref{19jun1208})
\begin{equation}\label{20jan1920}
\begin{aligned}
x_\eps^\lambda&=\log\eps^{-1}+\log\lambda^{-1}-\log(\log\eps^{-1}+\log\lambda^{-1})
-\log \bar c- 
m_0\farc{\log(\log\eps^{-1}+\log\lambda^{-1})}{\log\eps^{-1}+\log\lambda^{-1}}\\
&-\frac{m_1}{\log\eps^{-1}+\log\lambda^{-1}}-
\frac{m_2\log\bar c}{\log\eps^{-1}+\log\lambda^{-1}}
- \farc{m_3\bar c_1}{\bar c(\log\eps^{-1}+\log\lambda^{-1})}
+O(\ell_\eps^{-1-\delta})\\
&=x_\eps-\log\lambda-\farc{\log\lambda^{-1}}{\log\eps^{-1}}+O(\ell_\eps^{-1-\delta}).
\end{aligned}
\end{equation} 
Comparing (\ref{20jan1918}) and (\ref{20jan1920}) we see that we should have
\begin{equation}\label{20jan1922}
m_2=-1,
\end{equation}
that, in view of (\ref{20jan1916}), implies that $m_3=1$. 

\subsubsection*{Some preliminary transformations and the self-similar variables}

The conclusion of Proposition~\ref{prop-jun2502} follows from a series of 
changes of variables that we now describe. 
We  first  go into the moving frame, writing solution to (\ref{jul610})-(\ref{jul612})
as
\begin{equation}\label{jun2008}
u_\eps(t,x)=\tilde u_\eps(t,x-2t+\farc{3}{2}\log (t+1)).
\end{equation}
The function $\tilde u_\eps(t,x)$ satisfies
\begin{eqnarray}\label{jun2060}
\pdr{\tilde u_\eps}{t}-\Big(2-\farc{3}{2(t+1)}\Big)\pdr{\tilde u_\eps}{x}=\pdrr{\tilde u_\eps}{x}+ 
\tilde u_\eps- \tilde u_\eps^2.
\end{eqnarray}
Next, we take out the exponential decay factor, writing
\begin{equation}\label{jun2062}
\tilde u_\eps(t,x)=e^{-x}z_\eps(t,x),
\end{equation}
which gives
\begin{eqnarray}\label{jun2064}
\pdr{z_\eps}{t}-\farc{3}{2(t+1)}\Big(z_\eps-\pdr{z_\eps}{x}\Big)=\pdrr{z_\eps}{x} 
-e^{-x}z_\eps^2.
\end{eqnarray}
As (\ref{jun2064}) is a perturbation of the standard heat equation, it is helpful to pass to the self-similar variables:
\begin{equation}\label{jun2066}
z_\eps(t,x)=v_\eps(\log (t+1), \farc{x}{\sqrt{t+1}}).
\end{equation}
The function $v_\eps(\tau,\eta)$ is the solution of
\begin{eqnarray}\label{jun2068}
\pdr{v_\eps}{\tau}-\frac{\eta}{2}\pdr{v_\eps}{\eta}-\pdrr{v_\eps}{\eta}-
\farc{3}{2}v_\eps+\farc{3}{2}e^{-\tau/2}\pdr{v_\eps}{\eta}
+e^{\tau-\eta\exp(\tau/2)}v_\eps^2=0,
\end{eqnarray}
with the initial condition 
\begin{equation}\label{jul602}
v_\eps(0,\eta)=\eps e^\eta\phi_0(\eta).
\end{equation}
 In order to get rid of the pre-factor $\eps$ in the initial condition (\ref{jul602}), 
and also to adjust the zero-order term in (\ref{jun2068}), 
it is convenient to represent $v_\eps(\tau,\eta)$ as
\begin{equation}\label{jun2070}
v_\eps(\tau,\eta)=\eps v_1(\tau,\eta)e^{\tau/2}.
\end{equation}
Here, $v_1(\tau,\eta)$ is the solution of
\begin{eqnarray}\label{jun2072}
\pdr{v_1}{\tau}-\frac{\eta}{2}\pdr{v_1}{\eta}-\pdrr{v_1}{\eta}- 
v_1+\farc{3}{2}e^{-\tau/2}\pdr{v_1}{\eta}
+\eps e^{3\tau/2-\eta\exp(\tau/2)}v_1^2=0,
\end{eqnarray}
with the initial condition
\begin{equation}\label{jun2080}
v_1(0,\eta)=e^{\eta}\phi_0(\eta).
\end{equation}
 %
%

The next, and last, in this chain of preliminary transformations 
is to eliminate the pre-factor~$\eps$ in the last term in (\ref{jun2072}). We choose 
\begin{equation}\label{jun2090}
\beta(\tau)=e^{-\tau/2}\log\eps,
\end{equation}
so that
\begin{equation}\label{jun2088}
\eps e^{3\tau/2-\eta\exp(\tau/2)}= e^{3\tau/2-(\eta-\beta(\tau))\exp(\tau/2)},
\end{equation}
and make a change of the spatial variable:
\begin{equation}\label{jun2092}
v_1(\tau,\eta)=r_\eps(\tau,\eta-\beta(\tau)).
\end{equation}
The function $r_\eps$ satisfies:
\begin{eqnarray}\label{jun20104bis}
\pdr{r_\eps}{\tau}-\frac{\eta}{2}\pdr{r_\eps}{\eta}
-\pdrr{r_\eps}{\eta}- r_\eps
+\farc{3}{2}e^{-\tau/2}\pdr{r_\eps}{\eta}
+ e^{3\tau/2-\eta\exp(\tau/2)}r_\eps^2=0,
\end{eqnarray}
with the initial 
condition 
\begin{equation}\label{19jul2504}
r_\eps(0,\eta)=\psi_0(\eta-\ell_\eps),
\end{equation}
with $\ell_\eps$ as in (\ref{19jun1210}), and
\begin{equation}\label{19jul2502}
\psi_0(\eta)=e^\eta\phi_0(\eta).
\end{equation}
This, with a slight abuse of notation,
is exactly (\ref{jun20104}).
Note that $r_\eps$ depends on $\eps$ only through $\ell_\eps$ as it appears in
the initial condition. We will interchangeably, with some abuse of notation use $r_\eps(t,x)$
and $r_{l_\eps}(t,x)$ for the same object. 
 
As far as the asymptotics of $r_\eps(\tau,\eta)$ and its connection to the Bramson shift are 
concerned, 
it was shown in \cite{NRR1} that there exists a constant~$v_\infty(\eps)>0$
so that the solution $v_\eps(\tau,\eta)$ of (\ref{jun2068})
has the asymptotics
\begin{equation}\label{may1002}
v_\eps(\tau,\eta)\sim v_\infty(\eps)\eta e^{-\eta^2/4}e^{\tau/2},
\hbox{ as $\tau\to+\infty$, for $\eta>0$}.
\end{equation}
and the Bramson shift is given by
\begin{equation}\label{may1004}
x_\eps=-\log v_\infty(\eps).
\end{equation}
The corrseponding long-time asymptotics for the function $v_1(\tau,\eta)$,
the solution to (\ref{jun2072}) is
\begin{equation}\label{may1006}
v_1(\tau,\eta)\sim \tilde v_\infty(\eps)\eta e^{-\eta^2/4},
\hbox{ as $\tau\to+\infty$, for $\eta>0$},
\end{equation}
with
\begin{equation}\label{may1008}
\tilde v_\infty(\eps) =\eps v_\infty(\eps),
\end{equation}
and the asymptotics for $r_\eps$ is
\begin{equation}\label{may1012bis}
r_\eps(\tau,\eta)=v_1(\tau,\eta+\beta(\tau))\sim \tilde v_\infty(\eps)(\eta+\beta(\tau)) 
e^{-(\eta+\beta(\tau)^2/4}\sim
\tilde v_\infty(\eps)\eta e^{-\eta^2/4},\hbox{ as $\tau\to+\infty$, for $\eta>0$},
\end{equation}
so that
\begin{equation}\label{may1014}
r_\infty(\ell_\eps)=\tilde v_\infty(\eps)=\eps v_\infty(\eps),
\end{equation}
and the Bramson shift is
\begin{equation}\label{jun2510}
x_\eps=-\log v_\infty(\eps)=
\log\eps^{-1}-\log r_\infty(\ell_\eps).
\end{equation} 
This finishes the proof of Proposition~\ref{prop-jun2502}.~$\Box$

\section{Connection to the linear Dirichlet problem}\label{sec:linear}

Before giving the proof of Proposition~\ref{prop-jun2504},
let us recall the intuition that leads to
the long-time asymptotics~(\ref{may1012}) 
for the solution of (\ref{jun20104}), and also explain how the asymptotics (\ref{jun2502}) comes about. 
The key point is that we may think of (\ref{jun20104}) as a linear equation with the factor
\begin{equation}\label{nov802}
e^{3\tau/2-\eta\exp(\tau/2)}r_\ell(\tau,\eta)
\end{equation}
in the last term in its right side playing the role of an absorption coefficient. 
Disregarding our lack of information about $r_\ell(\tau,\eta)$ that enters (\ref{nov802}),
we expect that
when $\tau\gg 1$  this term is "extremely large" for $\eta<0$
and "extremely small" for $\eta>0$. Thinking again of~(\ref{jun20104}) as a linear
equation for~$r_\ell(\tau,\eta)$, the former means that~$r_\ell(\tau,\eta)$ 
is very small for $\eta<0$,
while the latter indicates that~$r_\ell(\tau,\eta)$ essentially solves a linear problem
for $\eta>0$. The drift term in (\ref{jun20104}) with the pre-factor $e^{-\tau/2}$
is also very small at large times. 
Thus, if we take some $T\gg 1$, then for $\tau\ge T$, a good approximation
to~(\ref{jun20104}) is the linear Dirichlet problem
\begin{equation}\label{may2414}
\begin{aligned}
&\pdr{\zeta_\ell}{\tau}-\frac{\eta}{2}\pdr{\zeta_\ell}{\eta}
-\pdrr{\zeta_\ell}{\eta}- \zeta_\ell=0,~~\tau>T,~\eta>0\\
&\zeta_\ell(\tau,0)=0,\\
&\zeta_\ell(T,\eta)=r_\ell(T,\eta).
\end{aligned}
\end{equation}
In other words, one would solve the full nonlinear problem on the whole line
only until a large time~$T\gg 1$, and for $\tau>T$ simply solve the linear Dirichlet problem
(\ref{may2414}). 
It is easy to see that 
\begin{equation}\label{may2416}
\bar\zeta(\eta)=\eta e^{-\eta^2/4},
\end{equation}
is a steady solution to (\ref{may2414}). In addition, the operator  
\begin{equation}\label{may2420}
{\cal L}u=\frac{\partial^2u}{\partial\eta^2}+\frac{\eta}{2}\pdr{u}{\eta}+u,~~\eta>0,
\end{equation}
with the Dirichlet boundary condition at $\eta=0$ has a discrete spectrum. It follows
that $\zeta_\ell(\tau,\eta)$ has the long time asymptotics
\begin{equation}\label{may2418}
\zeta_\ell(\tau,\eta)\sim\zeta_\infty(\ell)\eta e^{-\eta^2/4},~~\hbox{$\tau\to+\infty$}.
\end{equation}
As the integral
\begin{equation}\label{may2422}
\int_0^\infty\eta\zeta_\ell(\tau,\eta)d\eta=\int_0^\infty\eta\zeta_\ell(T,\eta)d\eta
\end{equation}
is conserved, the coefficient $\zeta_\infty(\ell)$ is determined by the relation
\begin{equation}\label{may2424}
\zeta_\infty(\ell)\int_0^\infty \eta^2 e^{-\eta^2/4}d\eta=
\int_0^\infty\eta\zeta_\ell(T,\eta)d\eta,
\end{equation}
so that
\begin{equation}\label{may2426}
\zeta_\infty(\ell)=\farc{1}{\sqrt{4\pi}}\int_0^\infty\eta\zeta_\ell(T,\eta)d\eta=
\farc{1}{\sqrt{4\pi}}\int_0^\infty\eta r_\ell(T,\eta)d\eta.
\end{equation}
As we expect $\zeta_\ell(\tau,\eta)$ and $r_\ell(\tau,\eta)$ to be close if $T$ is 
sufficiently large, we should have an approximation
\begin{equation}\label{may2502}
\zeta_\infty(\ell)\approx r_\infty(\ell),
\end{equation}
if $T\gg 1$. This, in turn, implies that
\begin{equation}\label{may2504}
r_\infty(\ell)=\lim_{\tau\to+\infty}
\farc{1}{\sqrt{4\pi}}\int_0^\infty\eta r_\ell(\tau,\eta)d\eta.
\end{equation}
This informal argument is made rigorous in~\cite{NRR1}. 

The limit in the right side of (\ref{may2504})
is an implicit functional of the initial
conditions for the nonlinear problem (\ref{jun20104}), and the evolution of the solution
in the initial time layer, before the linear approximation kicks in, is difficult to control,
so that there is no explicit expression for $r_\infty(\ell)$.
In the present setting, however, the initial condition $r_\ell(0,\eta)$ in (\ref{jun20104})
is shifted to the right by $\ell\gg 1$. Therefore, at small times 
the solution is concentrated at $\eta\gg 1$, a region where the factor in front of the
nonlinear term in~(\ref{jun20104})
\begin{equation}\label{may2506}
\exp\Big(\frac{3\tau}{2}-\eta e^{\tau/2}\Big)\ll 1
\end{equation}
is very small even for $\tau=O(1)$. Hence, solutions to
the nonlinear equation (\ref{jun20104}) with the initial conditions (\ref{19jul2504})
should be well approximated, to the leading order, 
by the linear problem 
\begin{eqnarray}\label{may2508}
\pdr{\tilde r_\ell}{\tau}-\frac{\eta}{2}\pdr{\tilde r_\ell}{\eta}
-\pdrr{\tilde r_\ell}{\eta}- \tilde r_\ell
+\farc{3}{2}e^{-\tau/2}\pdr{\tilde r_\ell}{\eta}=0,~~\tilde r_\ell(0,\eta)=r_\ell(0,\eta),
\end{eqnarray}
even for small times. However, the solution "does not yet know" for "small" $\tau$ that there
is a large dissipative term in the nonlinear equation, or the Dirichlet boundary condition
in the linear version, and evolves "as if (\ref{may2508}) is posed for $\eta\in\Rm$". This
leads to exponential growth in $\tau$ until the solution spreads sufficiently far to the
left, close to $\eta=0$ and "discovers" the Dirichlet 
boundary condition (or the nonlinearity in the full nonlinear version). 
During
this "short time" evolution we have 
\begin{equation}\label{may2510}
\farc{d}{d\tau}\int\eta\tilde r_\ell(\tau,\eta)d\eta=
\farc{3}{2}e^{-\tau/2}\int\tilde r_\ell(\tau,\eta)d\eta.
\end{equation}
Unlike the first moment, the total mass in the right side does not grow as $\ell$ 
gets larger -- the shift of the initial condition to the right increases 
the first moment but not
the mass. Thus, the first moment of $r_\ell(\tau,\eta)$
will only change by a factor that is $o(1)$ during the "short time" evolution,
so that it is conserved to the leading order in $\ell$. 
The "long time" evolution following this initial time layer
is well 
approximated by the linear Dirichlet problem (\ref{may2414}) that
preserves the first moment. Thus, altogether, 
the first moment will not change to the leading order
if $\ell\gg 1$ is large, so that
\begin{equation}\label{may2514}
\lim_{\tau\to+\infty}\int_0^\infty \eta r_\ell(\tau,\eta)d\eta=(1+o(1))
\int_0^\infty \eta r_\ell(0,\eta)d\eta,~~\hbox{ as $\eps\to 0$},
\end{equation}
which leads to the explicit expression for $r_\infty(\ell)$ in terms
of the initial first moment:
\begin{equation}\label{may2512}
\begin{aligned}
r_\infty(\ell)&=(1+o(1))\frac{1}{\sqrt{4\pi}}\int_0^\infty \eta r_\ell(0,\eta)d\eta=
(1+o(1))\frac{1}{\sqrt{4\pi}}\int_0^\infty \eta e^{\eta-\ell}
\phi_0(\eta-\ell)d\eta\\
&=(1+o(1))\frac{1}{\sqrt{4\pi}}\int_{-\ell}^\infty (\eta+\ell) 
e^{\eta}\phi_0(\eta)d\eta=(1+o(1))\frac{\ell}{\sqrt{4\pi}}
\int_{-\infty }^\infty  e^{\eta}\phi_0(\eta)d\eta\\
&=\bar c(1+o(1)) \ell,
\end{aligned}
\end{equation}
which is (\ref{jun2502}). This very informal argument is behind the reason why we can 
describe the Bramson shift so explicitly for $\eps\ll 1$, which corresponds to $\ell\gg 1$.
The rest of the proof of Proposition~\ref{prop-jun2504} formalizes
this argument by providing matching upper and lower bounds on the limit in the right
side of (\ref{may2504}).  
 
\subsection*{An approximate solution to the adjoint linear problem}

In order to improve on the approximate conservation law (\ref{may2510}) let us make
the following observation. Let us  set 
\begin{equation}\label{19jul2602}
Q_k(\tau,\eta)=\eta+k\bar\psi(\eta)e^{-\tau/2},
\end{equation}
with
\begin{equation}\label{jun2812}
\bar\psi(\eta)=\int_0^\eta e^{z^2/4}\int_{z}^\infty e^{-y^2/4}dy dz,
\end{equation}
and consider a solution to the linear Dirichlet problem
\begin{eqnarray}\label{jul192604}
\pdr{p}{\tau}-\frac{\eta}{2}\pdr{p}{\eta}
-\pdrr{p}{\eta}- p
+ke^{-\tau/2}\pdr{p}{\eta}=f(\tau,\eta),~~p(\tau,0)=0.
\end{eqnarray}
\begin{lemma}\label{lem-jul2602}
We have
\begin{equation}\label{19jul2606}
\farc{d}{d\tau}\int_0^\infty Q_k(\tau,\eta)p(\tau,\eta)d\eta= k^2e^{-\tau}\int_0^\infty
\pdr{\bar\psi(\eta)}{\eta}p(\tau,\eta)d\eta+\int_0^\infty f(\tau,\eta) 
Q_k(\tau,\eta)d\eta.
\end{equation}
\end{lemma}
{\bf Proof.} Note that
\begin{equation}\label{19jul2612}
\begin{aligned}
\frac{d}{d\tau}&\int_0^\infty \! Q_k(\tau,\eta)p(\tau,\eta)d\eta=
\int_0^\infty\Big(\pdr{Q_k }{\tau}p  + 
Q_k \Big[
\frac{\eta}{2}\pdr{p}{\eta}
+\pdrr{p}{\eta}+p
-ke^{-\tau/2}\pdr{p}{\eta}+f\Big]\Big) d\eta\\
&=\int_0^\infty\Big( p\Big[\pdr{Q_k}{\tau}+
\pdrr{Q_k}{\eta}-\pdr{}{\eta}\Big(\farc{\eta}{2} Q_k\Big)
+Q_k 
+ke^{-\tau/2}\pdr{Q_k }{\eta}\Big]+f Q_k \Big) d\eta.
\end{aligned}
\end{equation}
It is easy to check that the function $\bar\psi(\eta)$ 
is a solution to
\begin{equation}\label{19jul2610}
\begin{aligned}
&\farc{\eta}{2}\pdr{\bar\psi}{\eta}  
-\pdrr{\bar\psi}{\eta} -1=0,~~\bar\psi(0)=0.
\end{aligned}
\end{equation}
With this, we can compute that the function $Q_k(\tau,\eta)$ 
satisfies
\begin{equation}\label{19jul2512}
\begin{aligned}
&\pdr{Q_k}{\tau}-\pdr{}{\eta}\Big(\frac{\eta}{2}Q_k\Big)
+\pdrr{Q_k}{\eta}+ Q_k
+ke^{-\tau/2}\pdr{Q_k}{\eta}= ke^{-\tau/2}
-\farc{k}{2}e^{-\tau/2}\bar\psi(\eta)-\farc{k}{2}\bar\psi(\eta)e^{-\tau/2}\\
&-
\farc{k\eta}{2}\pdr{\bar\psi(\eta)}{\eta}e^{-\tau/2}
+k\pdrr{\bar\psi(\eta)}{\eta}e^{-\tau/2}+k\bar\psi(\eta) e^{-\tau/2}
+k^2 \pdr{\bar\psi(\eta)}{\eta}e^{-\tau}\\
&=k\Big(1-\farc{\eta}{2}\pdr{\bar\psi(\eta)}{\eta}+\pdrr{\bar\psi(\eta)}{\eta} \Big)
e^{-\tau/2}+k^2\pdr{\bar\psi(\eta)}{\eta}e^{-\tau}=k^2\pdr{\bar\psi(\eta)}{\eta}e^{-\tau}.
\end{aligned}
\end{equation}
Using this in (\ref{19jul2612}) gives (\ref{19jul2606}).~$\Box$

Note that for $0\le\eta\le 1$ we have
\begin{equation}\label{jun2816}
\bar\psi(\eta)\le C\eta.
\end{equation}
For the asymptotics of $\bar\psi(\eta)$ for large $\eta\gg 1$, note that
\begin{equation}\label{19jul2614}
\int_z^\infty e^{-y^2/4}dy=-\int_z^\infty\farc{2}{y}\Big(e^{-y^2/4}\Big)'dy=
\frac{2}{z}e^{-z^2/4}-\int_z^\infty \farc{2}{y^2}e^{-y^2/4}dy.
\end{equation}
It follows that for $\eta\ge 1$ we have
\begin{equation}\label{19jul26114}
\begin{aligned}
\bar\psi(\eta)&=\bar\psi(1)+\int_1^\eta e^{z^2/4}\int_{z}^\infty e^{-y^2/4}dy dz
=\bar\psi(1)+2\log \eta-2\int_1^\eta e^{z^2/4}\int_z^\infty \farc{1}{y^2}e^{-y^2/4}dy\\
&=2\log \eta+g(\eta),~~\eta\ge 1,
\end{aligned}
\end{equation}
where $g(\eta)$ is a bounded function such that
\begin{equation}\label{nov804}
g_\infty=\lim_{\eta+\infty} g(\eta)=
\int_0^1 e^{z^2/4}\int_z^\infty e^{-y^2/4}dy dz
-2\int_1^\infty e^{z^2/4}\int_z^\infty \farc{1}{y^2}e^{-y^2/4}dy.
\end{equation}
This is how the constant $g_\infty$ appears in (\ref{20apr1440}) and
in Theorem~\ref{conj-jun1102-19}. 

%

\section{The proof of Proposition~\ref{prop-jun2504}}\label{sec:prop-lln}

\subsection{An upper bound for the first moment} 
 
In this section, we prove an upper bound for  $r_\infty(\ell)$. 
\begin{lemma}\label{lem-may2502}
There exists $K>0$ so that we have
\begin{equation}\label{may2520}
r_\infty(\ell)=\limsup_{\tau\to+\infty}
\farc{1}{\sqrt{4\pi}}\int_0^\infty\eta r_\ell(\tau,\eta)d\eta\le 
\bar c\ell +K\log\ell
~~\hbox{for $\ell\ge 2$},
\end{equation}
with $\bar c$ as in  (\ref{jul620}).
\end{lemma}
 
\subsubsection*{Reduction to a Dirichlet problem} 

We first bound the solution to (\ref{jun20104}) by a solution to the
linear Dirichlet problem, up to a relatively small error. We start with two observations.
First, the solution of the original KPP problem (\ref{jul610}) satisfies
$u(t,x)\le 1$, hence the function~$v_\eps(\tau,\eta)$ defined in
(\ref{jun2066}) satisfies 
\[
v_\eps(\tau,\eta)\le e^{\eta e^{\tau/2}}.
\]
Retracing our changes of variables, we deduce that
 $v_1(\tau,\eta)$ defined in (\ref{jun2070}) satisfies
\begin{equation}\label{jul1002}
v_1(\tau,\eta)= \eps^{-1} v_\eps(\tau,\eta) e^{-\tau/2}
\le \eps^{-1}e^{-\tau/2+\eta e^{\tau/2}},
\end{equation}
and $r_\eps(t,x)$ introduced in (\ref{jun2092}) obeys
\begin{equation}\label{jul1004}
r_{\ell_\eps}(\tau,\eta)=v_1(\tau,\eta+\beta(\tau))\le
\eps^{-1}e^{(\eta-e^{-\tau/2}\log\eps^{-1}) e^{\tau/2}}e^{-\tau/2}=
e^{\eta e^{\tau/2}-\tau/2}.
\end{equation}
It follows that at the boundary $\eta=0$ we have 
\begin{equation}\label{jul1006}
0<r_\ell(\tau,0)\le e^{-\tau/2},~~\hbox{ for all $\tau>0$,}
\end{equation}
so that for $\tau\gg 1$ the function $r_\ell$ does satisfy an approximate Dirichlet boundary
condition at~$\eta=0$. 
However, the bound (\ref{jul1006}) is very poor for $\tau=O(1)$ -- recall that 
the initial condition is located at distance~$\ell\gg 1$ away from the
origin, so the solution remains small near $\eta=0$ for some time $\tau\gg 1$. 
In particular, as a first step,
we can bound $r_\ell(\tau,\eta)$ from above by
the solution to the linear problem on the whole line:
\begin{eqnarray}\label{jul1202}
&&\pdr{\bar R_\ell}{\tau} -
\frac{\eta}{2}\pdr{\bar R_\ell}{\eta}-
\pdrr{\bar R_\ell}{\eta}- \bar R_\ell
+\farc{3}{2}e^{-\tau/2}\pdr{\bar R_\ell}{\eta}=0,~~\eta\in\Rm,
\nonumber\\
&&\bar R_\ell(0,\eta)=r_{\ell}(0,\eta). 
\end{eqnarray}
A change of variables 
\[
\bar R_\ell(\tau,\eta)=Q_\ell(\tau,\eta-\farc{3}{2}\tau e^{-\tau/2})
\]
leads to the standard heat equation in the self-similar variables
\begin{eqnarray}\label{jul1202bis}
&&\pdr{Q_\ell}{\tau} -
\frac{\eta}{2}\pdr{Q_\ell}{\eta}-
\pdrr{Q_\ell}{\eta}- Q_\ell
=0,~~\eta\in\Rm,
\nonumber\\
&&Q_\ell(0,\eta)=r_{\ell}(0,\eta). 
\end{eqnarray}
Thus, the function $\bar R_\ell(\tau,\eta)$ can be written explicitly as
\begin{equation}\label{jul1204}
\begin{aligned}
\bar R_{\ell}(\tau,\eta)&=e^\tau\int G(e^{\tau}-1,(\eta-\farc{3}{2}\tau e^{-\tau/2}) e^{\tau/2}
-y)r_\ell(0,y)dy\\
&=e^\tau\int G(e^{\tau}-1,\eta e^{\tau/2}-\farc{3}{2}\tau  -y)r_\ell(0,y)dy.
\end{aligned}
\end{equation}
Here, $G(t,x)$ is the standard heat kernel:
\begin{equation}\label{nov806}
G(t,x)=\frac{1}{\sqrt{4\pi t}}e^{-|x|^2/(4t)}.
\end{equation}
As $r_\ell(0,\eta)=\psi_0(\eta-\ell)$, and $\psi_0(\eta)$ satisfies $\psi_0(\eta)=0$
for $\eta>L_0$ and $\psi_0(\eta)\le C e^\eta$ for $\eta<0$,
we have
\begin{equation}\label{jul1206}
\begin{aligned}
\bar R_{\ell}(\tau,0)&= e^\tau\int G(e^{\tau}-1,y+\farc{3}{2}\tau)r_\ell(0,y)dy=
e^\tau\int G(e^{\tau}-1,y+\ell+\farc{3}{2}\tau)\psi_0(y)dy\\
&\le Ce^\tau\int_{-\infty}^{L_0} G(e^{\tau}-1,y+\ell+\farc{3}{2}\tau)e^ydy
= C e^\tau e^{-\ell-(3/2)\tau}\int_{-\infty}^{L_0+\ell+(3\tau/2)} G(e^{\tau}-1,y)e^ydy.
\end{aligned}
\end{equation}
Note that for any $L\in\Rm$ we have
\begin{equation}\label{nov1202}
\begin{aligned}
\int_{-\infty}^L G(t,y)e^ydy&=\farc{1}{\sqrt{4\pi t}}\int_{-\infty} ^{L}\exp\Big(-\farc{y^2}{4t}+y\Big)dy=
\farc{1}{\sqrt{4\pi }}\int_{-\infty} ^{L/\sqrt{t}}\exp\Big(-\farc{y^2}{4}+y\sqrt{t}-t+t\Big)dy\\
&=\farc{e^{t}}{\sqrt{4\pi }}\int_{-\infty} ^{L/\sqrt{t}}\exp\Big(-\Big(\farc{y}{2}-\sqrt{t}\Big)^2 \Big)dy
=\farc{e^{t}}{\sqrt{4\pi }}\int_{-\infty} ^{L/\sqrt{t}-2\sqrt{t}}\exp\Big(-\farc{y^2}{4} \Big) dy.
\end{aligned}
\end{equation}
We take $\delta>0$ sufficiently small, and 
consider two cases. First, if 
\[
\tau<\tau_1=\log\ell-\log(2-\delta),
\]
then we have, from (\ref{jul1206}) and (\ref{nov1202}), for  for $0<\tau<\tau_1$:
\begin{equation}\label{nov1204}
\begin{aligned}
\bar R_{\ell}(\tau,0)&\le Ce^{-\ell}e^{-\tau/2}e^{e^\tau-1} \le Ce^{-\ell}e^{-\tau_1/2}e^{e^{\tau_1}-1} 
=  Ce^{-\ell}e^{-\tau_1/2} e^{{\ell}}= C\ell^{-1/2}\le C\ell^{-1/4}e^{-\tau/4}. 
\end{aligned}
\end{equation}
On the other hand, if $\tau>\tau_1$,  
then, taking 
\begin{equation}\label{nov1218}
t=e^\tau-1,~~L=L_0+\ell+\farc{3\tau}{2},
\end{equation}
in (\ref{nov1202}), we see that for $\ell$ sufficiently large, the upper limit of integration
\begin{equation}\label{nov1220}
\begin{aligned}
\farc{L}{\sqrt{t}}-2\sqrt{t}&=\farc{1}{\sqrt{e^{\tau}-1}}
\Big(L_0+\ell+\farc{3\tau}{2}-2e^\tau+2\Big)
=\farc{1}{\sqrt{e^{\tau}-1}}\Big(L_0+(2-\delta)e^{\tau_1}+\farc{3\tau}{2}-2e^\tau+2\Big)\\
&\le -\frac{\delta}{2}e^{\tau/2}
\end{aligned}
\end{equation}
is very negative. In particular, the integral in the right side of (\ref{nov1202}) can be estimated as
\begin{equation}\label{nov1206}
\begin{aligned}
&\int_{-\infty} ^{L/\sqrt{t}-2\sqrt{t}}\exp\Big(-\farc{y^2}{4} \Big) dy\le C\Big|\farc{L}{\sqrt{t}}-2\sqrt{t}\Big|^{-1}
\exp\Big(-\farc{1}{4}\Big(\farc{L}{\sqrt{t}}-2\sqrt{t}\Big)^{2}\Big)\\
&=
 C\Big|\farc{L}{\sqrt{t}}-2\sqrt{t}\Big|^{-1}e^{-L^2/(4t)}e^Le^{-t}.
\end{aligned}
\end{equation}
Then, we have from (\ref{jul1206}), (\ref{nov1202}), (\ref{nov1220}) and (\ref{nov1206})
\begin{equation}\label{nov1208}
\begin{aligned}
&\bar R_{\ell}(\tau,0)\le  
Ce^{-\ell-\tau/2}e^t\Big|\farc{L}{\sqrt{t}}-2\sqrt{t}\Big|^{-1}e^{-L^2/(4t)}e^Le^{-t}\\
&\le C_\delta e^{-\ell-\tau/2}e^{-\tau/2}e^{L_0+\ell+3\tau/2}e^{-(L_0+\ell+3\tau/2)^2/(4(e^\tau-1))}
\le C_\delta e^{\tau/2}e^{-(L_0+\ell+3\tau/2)^2/(4(e^\tau-1))}\\
&\le C e^{-\tau/2}\le Ce^{-\tau/4}e^{-\tau_1/4}\le C\ell^{-1/4}e^{-\tau/4},   
\end{aligned}
\end{equation}
provided that $\tau\ge\tau_1$ and 
\begin{equation}\label{nov1210}
\tau\le (L_0+\ell+3\tau/2)^2/(4(e^\tau-1)).
\end{equation}
In particular, we can take
\begin{equation}\label{jul1208}
\begin{aligned}
\tau\le \tau_2&=2\log\ell-\log\log\ell-3,
\end{aligned}
\end{equation}
as long as $\ell$ is sufficiently large, because then we have 
\begin{equation}\label{nov1212}
4\tau(e^\tau-1)\le 4\tau e^\tau \le 2\log\ell\farc{4\ell^2}{e^3\log \ell} \le \ell^2,
\end{equation}
so that (\ref{nov1210}) holds. 
%
It follows that
\begin{equation}\label{jul1234}
r_\ell(\tau,0)\le C\ell^{-1/4}e^{-\tau/4},~~~ 0\le\tau\le\tau_2.
\end{equation}
Taking into account (\ref{jul1006}), we deduce that we also have
\begin{equation}\label{jul1232}
r_\ell(\tau,0)\le e^{-\tau/2}\le e^{-\tau_2/4}e^{-\tau/4}=
\farc{C(\log\ell)^{1/4}}{\ell^{1/2}}e^{-\tau/4}\le 
\farc{C }{\ell^{1/4}}e^{-\tau/4},~~\hbox{for $\tau>\tau_2$}.
\end{equation}
It follows that
the function $r_\ell(\tau,\eta)$, the solution to (\ref{jun20104}),
is bounded from above by the solution to the linear half-line problem
\begin{eqnarray}\label{jun714}
&&\pdr{\bar r_\ell}{\tau} -
\frac{\eta}{2}\pdr{\bar r_\ell}{\eta}-
\pdrr{\bar r_\ell}{\eta}- \bar r_\ell
+\farc{3}{2}e^{-\tau/2}\pdr{\bar r_\ell}{\eta}=0,~~\eta>0, 
\nonumber\\
&&\bar r_\ell(\tau,0)=C\ell^{-1/4} 
e^{-\tau/4},
\end{eqnarray}
with the initial condition  
\begin{equation}\label{jun716}
\bar r_\ell(0,\eta)=\psi_0(\eta-\ell).
\end{equation}
In order to deal with the small but non-zero boundary condition
in (\ref{jun714}), we make yet another change of variables:
\begin{equation}\label{jun718}
\bar r_\ell(\tau,\eta)=\bar q_\ell(\tau,\eta )+
C\ell^{-1/4} e^{-\tau/4} g(\eta ),
\end{equation}
with a smooth function $g(\eta)$ such that $g(0)=1$, 
$g(\eta)\ge 0$ and $g(\eta)=0$ for $\eta\ge 1$.
This leads to
\begin{eqnarray}\label{jun720bis}
&&\pdr{\bar q_\ell}{\tau}+
 \farc{3}{2} e^{-\tau/2} \pdr{\bar q_\ell}{\eta}
={\cal L}\bar q_\ell+
C\ell^{-1/4}G(\tau,\eta)e^{-\tau/4},~~\eta>0,
\nonumber\\
&&\bar q_\ell(\tau,0)=0,
\end{eqnarray}
with a uniformly 
bounded function $G(\tau,\eta)$ that is supported in $\eta\in[0,1]$ and is 
independent of~$\ell$. We recall that the operator ${\cal L}$ is defined in (\ref{may2420}). 
This change of variable does not affect the 
asymptotics of the first moment:
\begin{equation}\label{jun728}
\limsup_{\tau\to+\infty}\int_0^\infty \eta r_\ell(\tau,\eta)d\eta\le
\limsup_{\tau\to+\infty}\int_0^\infty \eta \bar r_\ell(\tau,\eta)d\eta=
\limsup_{\tau\to+\infty}\int_0^\infty \eta \bar q_\ell(\tau,\eta)d\eta.
\end{equation}
The initial condition for $\bar q_\ell$ is
\begin{equation}\label{jun726}
\bar q_\ell(0,\eta)=\bar r_\ell(0,\eta)-
C\ell^{-1/4}g(\eta).
\end{equation}
It is easy to see that
\begin{equation}\label{jul640}
\int_0^\infty \eta \bar q_\ell(\tau,\eta)d\eta\le 
\int_0^\infty \eta\bar p_\ell(\tau,\eta)d\eta+C_0\ell^{-1/4}. 
\end{equation}
Here, $C_0$ is a constant that depends neither on $\ell$ nor on $\phi_0$, 
and $\bar p_\ell(\tau,\eta)$ is the solution to the homogeneous problem
\begin{eqnarray}\label{jun720}
&&\pdr{\bar p_\ell}{\tau}+
 \farc{3}{2} e^{-\tau/2}\pdr{\bar p_\ell}{\eta}
={\cal L}\bar p_\ell,~~\eta>0,
\nonumber\\
&&\bar p_\ell(\tau,0)=0,\\
&&\bar p_\ell(0,\eta)=\psi_0(\eta-\ell).\nonumber
\end{eqnarray}
Let us now use Lemma~\ref{lem-jul2602}, with $k=3/2$: 
multiply (\ref{jun720}) by $Q_{3/2}(\tau,\eta)$ 
and integrate:
\begin{equation}\label{19jul2522}
\begin{aligned}
\frac{d}{d\tau}&\int_0^\infty Q_{3/2}(\tau,\eta) 
\bar p_\ell(\tau,\eta)d\eta
=\farc{9}{4}e^{-\tau}\int_0^\infty \bar p_\ell(\tau,\eta)\pdr{\bar\psi(\eta)}{\eta}d\eta
\le Ce^{-\tau}\int_0^\infty \frac{\bar p_\ell(\tau,\eta)}{1+\eta}d\eta\\
&\le Ce^{-\tau}\int_0^\infty  \bar p_\ell(\tau,\eta) d\eta.
\end{aligned}
\end{equation}
Note, however, that integrating (\ref{jun720}) gives a simple upper bound
\begin{equation}\label{19jul2524}
\begin{aligned}
&\frac{d}{d\tau}\int_0^\infty \bar p_\ell(\tau,\eta)d\eta=
-\pdr{\bar p_\ell(\tau,0)}{\eta}+\frac{1}{2}\int_0^\infty \bar p_\ell(\tau,\eta)d\eta
\le \frac{1}{2}\int_0^\infty \bar p_\ell(\tau,\eta)d\eta,
\end{aligned}
\end{equation}
implying a trivial and very poor upper bound 
\begin{equation}\label{19jul2526}
\begin{aligned}
& \int_0^\infty \bar p_\ell(\tau,\eta)d\eta\le C_0 e^{\tau/2}.
\end{aligned}
\end{equation}
We remind that we use the notation $C_0$ for various constants that do not depend on $\eps$.
Using this estimate in the right side of (\ref{19jul2522}) gives
\begin{equation}\label{19jul2528}
\begin{aligned}
&\int_0^\infty Q_{3/2}(\tau,\eta)  
\bar p_\ell(\tau,\eta)d\eta\le  \int_0^\infty Q_{3/2}(0,\eta)  
\bar p_\ell(0,\eta)d\eta+C_0.
\end{aligned}
\end{equation}
Recalling the definition (\ref{19jul2602}) of $Q_k(\tau,\eta)$ and passing to the limit
$\tau\to+\infty$ gives
%
\begin{equation}\label{may1042}
\begin{aligned}
\bar p_\infty(\ell)&=\lim_{\tau\to+\infty}\farc{1}{\sqrt{4\pi}}
\int_0^\infty\eta\bar p_\ell(\tau,\eta)d\eta
\le \farc{1}{\sqrt{4\pi}}\int_0^\infty [\eta+\farc{3}{2}\bar\psi(\eta)] \psi_0(\eta-\ell)d\eta
+C_0\\
&=
\farc{1}{\sqrt{4\pi}}\int_{-\ell}^\infty\eta e^\eta\phi_0(\eta)d\eta+\farc{\ell}{\sqrt{4\pi}}
\int_{-\ell}^\infty e^\eta \phi_0(\eta)d\eta+ \farc{3}{2 } 
\int_{-\ell}^\infty\bar\psi(\eta+\ell)e^\eta\phi_0(\eta)\farc{d\eta}{\sqrt{4\pi}}+C_0\\
&\le \bar c\ell+ {3\bar c} \log\ell+C_1 
,\hbox{ for $\ell\ge 2$},
\end{aligned}
\end{equation}
with $\bar c$ as in (\ref{jul620}), and a constant $C_1$ that may depend on $\phi_0$.
We used (\ref{19jul26114}) in the last inequality above.
The conclusion of Lemma~\ref{lem-may2502} now follows. 

\subsection{A lower bound for  the first moment}

We now prove a lower bound for  the first moment matching the upper bound
in Lemma~\ref{lem-may2502}, to the leading order in $\ell$, which will finish the 
proof of the asymptotics (\ref{jun2502}) in Proposition~\ref{prop-jun2504}.  
\begin{lemma}\label{lem-may2504}
There exists a constant $K_1$ so that 
\begin{equation}\label{may2524}
\liminf_{\tau\to+\infty}
\farc{1}{\sqrt{4\pi}}\int_0^\infty\eta r_\ell(\tau,\eta)d\eta\ge 
\bar c\ell+K_1\log\ell,~~\hbox{for $\ell\ge 2 $.}
\end{equation}
\end{lemma}
Let us go back to (\ref{jun20104}):
\begin{eqnarray}\label{jun2110}
\pdr{r_\ell}{\tau}-\frac{\eta}{2}\pdr{r_\ell}{\eta}
-\pdrr{r_\ell}{\eta}- r_\ell
+\farc{3}{2}e^{-\tau/2}\pdr{r_\ell}{\eta}
+ e^{3\tau/2-\eta\exp(\tau/2)}r_\ell^2=0,
\end{eqnarray}
with the initial condition  
\begin{equation}\label{jun2112}
r_\ell(0,\eta)=e^{\eta-\ell }\phi_0(\eta-\ell ).
\end{equation}

We would like first to get rid of the nonlinearity in the last term in the left side of
(\ref{jun2110}). This is done as follows: recall that we have
\begin{equation}\label{jul1630}
r_\ell(\tau,\eta)\le \bar R(\tau,\eta),
\end{equation}
where $\bar R_\ell(\tau,\eta)$ is the solution to the linear whole line problem (\ref{jul1202}). Note that (\ref{jul1204})
implies that 
\begin{equation}\label{19jul2540}
\bar R_{\ell}(\tau,\eta)\le Ce^{\tau/2}\int \psi_0(y)dy= Ce^{\tau/2}.
\end{equation}
Thus, 
a lower bound for $r_\ell$ is the solution to
\begin{eqnarray}\label{jun2220}
\pdr{\tilde r_\ell}{\tau}-\frac{\eta}{2}\pdr{\tilde r_\ell}{\eta}
-\pdrr{\tilde r_\eps}{\ell}- \tilde r_\ell
+\farc{3}{2}e^{-\tau/2}\pdr{\tilde r_\ell}{\eta}
+ Ce^{2\tau-\eta\exp(\tau/2)}  \tilde r_\ell=0,~~~\eta> m\tau e^{-\tau/2},
\end{eqnarray}
with the boundary condition 
\begin{equation}\label{jun2224}
\tilde r_\ell(\tau, m\tau e^{\tau/2})=0, 
\end{equation}
and the initial condition 
\begin{equation}\label{jun2222}
\tilde r_\ell(0,\eta)= r_\ell(0,\eta),
\end{equation}
and with $m>2$.
Next, we shift, setting 
\[
\tilde p_\ell(\tau,\eta)=h_\ell(\tau,\eta-m\tau e^{-\tau/2}),
\]
so that $h_\ell(\tau,\eta)$ satisfies
\begin{equation}\label{jun2614}
\pdr{h_\ell}{\tau}
-\frac{\eta}{2} \pdr{h_\ell}{\eta}
-\pdrr{h_\ell}{\eta}- h_\ell
+\Big(\farc{3}{2}-m\Big)e^{-\tau/2}\pdr{h_\ell}{\eta}+ 
Ce^{(2-m)\tau -\eta e^{\tau/2}} 
h_\ell =0, 
\end{equation}
with the boundary condition 
\begin{equation}\label{jun2610}
h_\ell(\tau,0)=0,
\end{equation}
and the initial condition 
\begin{equation}\label{jun2612}
h_\ell(0,\eta)= r_\ell(0,\eta),~~\eta>0.
\end{equation}
We multiply (\ref{jun2614}) by 
\begin{equation}\label{19jul2542}
\tilde Q_m(\tau,\eta)=\eta+\Big(\farc{3}{2}-m\Big)\bar\psi(\eta)e^{-\tau/2},
\end{equation}
and integrate, using Lemma~\ref{lem-jul2602}: 
\begin{equation}\label{jun2616}
\begin{aligned}
\farc{d}{d\tau}&\int_0^\infty \tilde Q_m(\tau,\eta) h_\ell(\tau,\eta)d\eta=
-
C\int_0^\infty e^{ (2-m)\tau-\eta e^{\tau/2} }
Q_m(\tau,\eta) h_\ell(\tau,\eta) d\eta
\\
&+\Big(\farc{3}{2}-m\Big)^2
e^{-\tau}\int_0^\infty h_\ell(\tau,\eta)\pdr{\bar\psi(\eta)}{\eta}d\eta
\ge  
-Ce^{(3/2 -m)\tau} 
\int_0^\infty \eta e^{\tau/2}e^{-\eta e^{\tau/2}}
h_\ell(\tau,\eta) d\eta\\
&\ge 
-Ce^{(3/2-m)\tau  } 
\int_0^\infty 
h_\ell(\tau,\eta)d\eta.
\end{aligned}
\end{equation}
We used the fact that $m>2$ in the second inequality above, and that $\bar\psi$ is non-negative and increasing. 

Integrating (\ref{jun2614}) using positivity of $h(\tau,\eta)$ for $\eta>0$ gives 
\[
\frac{d}{d\tau}\int_0^\infty h_\ell(\tau,\eta)d\eta -\farc{1}{2}\int_0^\infty h_\ell(\tau,\eta)d\eta+\pdr{h(\tau,0)}{\eta}\le 0.
\]
As $\partial h(\tau,0)/\partial \eta>0$, this implies a trivial bound
\begin{equation}\label{may1050}
\int_0^\infty h_\ell(\tau,\eta)d\eta\le 
e^{\tau/2}\int_0^\infty h_\ell(0,\eta)d\eta=\bar c e^{\tau/2}.
\end{equation}
As $m>2$ in (\ref{jun2616}), we then obtain  
\begin{equation}\label{19jul2616}
\int_0^\infty \tilde Q_m(\tau,\eta) h_\ell(\tau,\eta)d\eta\ge 
\int_0^\infty \tilde Q_m(0,\eta) h_\ell(0,\eta)d\eta-C,
\end{equation}
so that 
\begin{equation}\label{19jul2618} 
\begin{aligned}
\bar p_\infty(\ell)&=\lim_{\tau\to+\infty}\farc{1}{\sqrt{4\pi}}
\int_0^\infty\eta\bar p_\ell(\tau,\eta)d\eta
\ge \liminf_{\tau\to+\infty}\farc{1}{\sqrt{4\pi}}
\int_0^\infty\eta  h_\ell(\tau,\eta)d\eta\\
&=\liminf_{\tau\to+\infty}\farc{1}{\sqrt{4\pi}}
\int_0^\infty \tilde Q_m(\tau,\eta)h_\ell(\tau,\eta)d\eta
\ge \farc{1}{\sqrt{4\pi}}\int_0^\infty [\eta+(\farc{3}{2}-m)\bar\psi(\eta)] 
\psi_0(\eta-\ell)d\eta
-C_0\\
&\ge  \bar c\ell-K \log\ell-C_1 
,\hbox{ for $\ell\ge 2$},
\end{aligned}
\end{equation}
with $\bar c$ as in (\ref{jul620}), and an appropriate $K\in\Rm$. 
This finishes the proof of the lower bound in
Lemma~\ref{lem-may2504}, and also completes the proof of estimate (\ref{jun2502}) in
Proposition~\ref{prop-jun2504}.~$\Box$

\section{The proof of Proposition~\ref{conj-jun1104-19}}\label{sec:5}
 
The proof of Proposition~\ref{conj-jun1104-19} is much more involved than that
of Proposition~\ref{prop-jun2504}. In this section, we outline the main
steps and state the auxiliary results needed in the proof.
 
We start with (\ref{jun20104bis}):  
\begin{eqnarray}\label{19jun1220}
\pdr{r_\ell}{\tau}-\frac{\eta}{2}\pdr{r_\ell}{\eta}
-\pdrr{r_\ell}{\eta}- r_\ell
+\farc{3}{2}e^{-\tau/2}\pdr{r_\ell}{\eta}
+ e^{3\tau/2-\eta\exp(\tau/2)}r_\ell^2=0,
\end{eqnarray}
with the initial condition  $r_\ell(0,\eta)=\psi_0(\eta-\ell)$. Recall that we are
interested in
\begin{equation}\label{19jun1222}
r_\infty(\ell)=\lim_{\tau\to+\infty}\farc{1}{\sqrt{4\pi}}
\int_0^\infty\eta r_\ell(\tau,\eta)d\eta=
\lim_{\tau\to+\infty}\farc{1}{\sqrt{4\pi}}
\int_0^\infty Q(\tau,\eta) r_\ell(\tau,\eta)d\eta.
\end{equation}
Here, $Q(\tau,\eta)=Q_{3/2}(\tau,\eta)$ is the approximate solution to the adjoint 
linear problem,  defined in~(\ref{19jul2602}) with~$k=3/2$.  
Multiplying (\ref{19jun1220}) by $Q(\tau,\eta)$
and integrating in $\eta$ gives, according to Lemma~\ref{lem-jul2602}:
\begin{equation}\label{19jul2630}
\farc{d}{d\tau}\int_0^\infty Q(\tau,\eta) r_\ell(\tau,\eta)d\eta=
-\int_0^\infty e^{3\tau/2-\eta\exp(\tau/2)}r_\ell^2(\tau,\eta)Q(\tau,\eta) d\eta+E_\ell(\tau),
\end{equation}
so that
\begin{equation}\label{19jun1226}
\begin{aligned}
r_\infty(\ell)&=Q_\ell+Y_\ell+{\cal E}_\ell,
\end{aligned}
\end{equation}
with
\be\label{20apr202}
\bal
Q_\ell&=\frac{1}{\sqrt{4\pi}}\int_0^\infty Q(0,\eta)r_\ell(0,\eta)d\eta,\\
Y_\ell&=-\frac{1}{\sqrt{4\pi}}\int_0^\infty\int_0^\infty e^{3\tau/2-\eta\exp(\tau/2)}
r_\ell^2(\tau,\eta)Q(\tau,\eta)
d\tau d\eta\\
{\cal E}_\ell&=\frac{1}{\sqrt{4\pi}}
\int_0^\infty E_\ell(\tau)d\tau.
\end{aligned}
\end{equation}
The error term
\begin{equation}\label{19aug2702}
E_\ell(\tau)=E_\ell^{(1)}(\tau)+E_\ell^{(2)}(\tau) 
\end{equation}
has two contributions. The first  
\begin{equation}\label{19jul2626}
E_\ell^{(1)}(\tau)=r_\ell(\tau,0)\pdr{Q(\tau,0)}{\eta}
\end{equation}
comes from the boundary at $\eta=0$ since we do not have 
$r_\ell(\tau,0)=0$ but only that $r_\ell(\tau,0)$ is small.
The second comes from the error term in Lemma~\ref{lem-jul2602},
because $Q(\tau,\eta)$ is only an approximate solution to the adjoint
problem and not an exact one:
\begin{equation}\label{19aug2704}
E_\ell^{(2)}(\tau)=\farc{9}{4}e^{-\tau}
\int_0^\infty\pdr{\bar\psi(\eta)}{\eta}r_\ell(\tau,\eta)d\eta\le Ce^{-\tau}\int_0^\infty\farc{r_\ell(\tau,\eta)}{1+\eta}d\eta.
\end{equation}

The linear term and the error term in (\ref{19jun1226}) are quite straightforward to evaluate and
estimate, respectively. The main difficulty will be in finding the precise asymptotics
of the nonlinear term in the right side of (\ref{19jun1226}).

\subsection*{The error term bound}

The error term in (\ref{19jun1226}) is bounded by the following lemma.
\begin{lemma}\label{lem-20apr202}
There exists $C>0$ so that
\be\label{20apr204}
{\cal E}_\infty(\ell)\le\farc{C}{\ell^{1/8}}.
\ee
\end{lemma}
{\bf Proof.} As we have seen in the proof
of Lemma~\ref{lem-may2502} -- see (\ref{jul1234}) and (\ref{jul1232}), 
we have an upper bound
\begin{equation}\label{19jul2630bis}
r_\ell(\tau,0)\le \bar r_\ell(\tau,0)\le \farc{C }{\ell^{1/4}}e^{-\tau/4}
\end{equation}
hence $E_\ell^{(1)}$ can be bounded as 
\begin{equation}\label{19jul2632}
\int_0^\infty E_\ell^{(1)}(\tau)d\tau\le \farc{C }{\ell^{1/4}}.
\end{equation}

We now estimate $E_\ell^{(2)}$. 
As in the proof of the upper bound in Lemma~\ref{lem-may2502}, see (\ref{jun718})-(\ref{jun720}),
we deduce from
(\ref{19jul2630bis}) that
\begin{equation}\label{19aug2918}
E_\ell^{(2)}(\tau)\le \farc{C }{\ell^{1/4}}e^{-\tau}
+Ce^{-\tau}\int_0^\infty\farc{\bar p_\ell(\tau,\eta)}{1+\eta}d\eta,
\end{equation}
where $\bar p_\ell(\tau,\eta)$ is the solution to (\ref{jun720}):
\begin{equation}\label{19aug2920}
\begin{aligned}
&\pdr{\bar p_\ell}{\tau}+
 \farc{3}{2} e^{-\tau/2}\pdr{\bar p_\ell}{\eta}
=\pdrr{\bar p_\ell}{\eta}+\frac{\eta}{2}\pdr{\bar p_\ell}{\eta}
 + \bar p_\ell,~~\eta>0,
\\
&\bar p_\ell(\tau,0)=0,\\
&\bar p_\ell(0,\eta)=\psi_0(\eta-\ell)=e^{\eta-\ell}\phi_0(\eta-\ell).
\end{aligned}
\end{equation}
The simple-minded bound (\ref{19jul2526}) implies that if we fix any $T>0$, then  
\begin{equation}\label{nov1302}
e^{-\tau}\int_0^\infty\farc{\bar p_\ell(\tau,\eta)}{1+\eta}d\eta\le C e^{-\tau/2}\le Ce^{-\tau/4}e^{-T/4}
\hbox{ for all $\tau>T$.}
\end{equation}
For short times $\tau<T$ we can take $L>0$ and use (\ref{19jul2526}) to write
\begin{equation}\label{nov1304}
\begin{aligned}
e^{-\tau}\int_0^\infty\farc{\bar p_\ell(\tau,\eta)}{1+\eta}d\eta &\le 
C e^{-\tau} \Big(\sup_{\eta\in[0,L]}\bar p_\ell(\tau,\eta)\Big)\log L+\farc{C}{L}e^{-\tau/2}\\
&\le C e^{-\tau} \Big(\sup_{\eta\in[0,L]}\bar R_\ell(\tau,\eta)\Big)\log L+\farc{C}{L}e^{-\tau/2}.
\end{aligned}
\end{equation}
Here, $\bar R_\ell(\tau,\eta)$ is the solution to the whole line problem (\ref{jul1202}):
\begin{eqnarray}\label{jul1202bis2}
&&\pdr{\bar R_\ell}{\tau} -
\frac{\eta}{2}\pdr{\bar R_\ell}{\eta}-
\pdrr{\bar R_\ell}{\eta}- \bar R_\ell
+\farc{3}{2}e^{-\tau/2}\pdr{\bar R_\ell}{\eta}=0,~~\eta\in\Rm,
\nonumber\\
&&\bar R_\ell(0,\eta)=r_{\ell}(0,\eta),
\end{eqnarray}
given explicitly by (\ref{jul1204}):
\begin{equation}\label{jul1204bis}
\begin{aligned}
\bar R_{\ell}(\tau,\eta)& =e^\tau\int G(e^{\tau}-1,\eta e^{\tau/2}-\farc{3}{2}\tau  -y)r_\ell(0,y)dy.
\end{aligned}
\end{equation}
In the proof of Lemma~\ref{lem-may2502} we only looked for the bound on $\bar R_\ell(\tau,\eta)$ at
$\eta=0$ but now we need to consider $\eta\in[0,L]$. Note that for $T$ and $L$ not too large, 
the function $\bar R(\tau,\eta)$ is increasing in $\eta$ for~$\eta\in[0,L]$ and $0<\tau<T$. Hence, we have
\begin{equation}\label{nov1306}
\begin{aligned}
\bar R_{\ell}(\tau,\eta)& \le \bar R_\ell(\tau,L)= e^\tau\int G(e^{\tau}-1,L e^{\tau/2}-\farc{3}{2}\tau  -y)r_\ell(0,y)dy,
~~\hbox{ for~$\eta\in[0,L]$ and $0<\tau<T$.}
\end{aligned}
\end{equation} 
As in (\ref{jul1206}) and (\ref{nov1204}), we get, as $\psi_0(y)=0$ for $y\ge L_0$:
\begin{equation}\label{nov1308}
\begin{aligned}
\bar R_{\ell}(\tau,L)&= e^\tau\int G(e^{\tau}-1,y+\farc{3}{2}\tau-L e^{\tau/2})r_\ell(0,y)dy=
e^\tau\int G(e^{\tau}-1,y+\ell+\farc{3}{2}\tau-L e^{\tau/2})\psi_0(y)dy\\
&\le Ce^\tau\int_{-\infty}^{L_0} G(e^{\tau}-1,y+\ell+\farc{3}{2}\tau-L e^{\tau/2})e^ydy\\
&= C e^\tau e^{-\ell-(3/2)\tau+Le^{\tau/2}}\int_{-\infty}^{L_0+\ell+(3\tau/2)-Le^{\tau/2}} G(e^{\tau}-1,y)e^ydy\\
&\le C  e^{-\ell-\tau/2+Le^{\tau/2}}e^{e^\tau-1}\le C  e^{-\ell-T/2+Le^{T/2}}e^{e^T-1}.
\end{aligned}
\end{equation}
Let us take $T=(1/2)\log\ell$ and $L=\ell^{1/2}$, which gives 
\begin{equation}\label{nov1310}
\begin{aligned}
\bar R_{\ell}(\tau,\eta)&\le   C  e^{-\ell-\log\ell/4+\ell^{3/4} +\ell^{1/2}-1} \le C e^{-\ell/2}.
\end{aligned}
\end{equation} 
Using this in (\ref{nov1304}), we obtain for $\tau<T$:
\begin{equation}\label{nov1312}
\begin{aligned}
&e^{-\tau}\int_0^\infty\farc{\bar p_\ell(\tau,\eta)}{1+\eta}d\eta \le
C e^{-\tau} e^{-\ell/2}\log\ell+\farc{C}{\ell^{1/2}}e^{-\tau/2}\hbox{ for all $\tau<T$},
\end{aligned}
\end{equation}
while (\ref{nov1302}) becomes 
\begin{equation}\label{nov1314}
e^{-\tau}\int_0^\infty\farc{\bar p_\ell(\tau,\eta)}{1+\eta}d\eta
\le Ce^{-\tau/4}\ell^{-1/8}
\hbox{ for all $\tau>T$.}
\end{equation}
It follows that
\begin{equation}\label{19aug2912}
\int_0^\infty E_\ell^{(2)}(\tau)d\tau\le\farc{C}{\ell^{1/8}} ,
\ee
and the conclusion of Lemma~\ref{lem-20apr202} follows.~$\Box$

\subsection*{The linear contribution}

The term $Q_\ell$ in (\ref{20apr202})
can be computed explicitly from the asymptotics for $Q(0,\eta)$ as $\eta\to \infty$ coming from~(\ref{19jul26114}):
\begin{equation}\label{19jun1236bis}
Q(0,\eta)=\eta+3\log\eta+\farc{3}{2}g_\infty+O\Big(\farc{1}{\eta}\Big)\hbox{ as $\eta\to+\infty$},
\end{equation}
with $g_\infty$ defined in (\ref{nov804}). This gives
\begin{equation}\label{jan806}
\begin{aligned}
Q_\ell=\farc{1}{\sqrt{4\pi}}\int_0^\infty(\eta+3\log\eta+\farc{3}{2}g_\infty)
\psi_0(\eta-\ell)d\eta+O\Big(\farc{1}{\ell}\Big)=\bar c\ell+3\bar c\log\ell
+\farc{3}{2}g_\infty\bar c+\bar c_1+
O(\ell^{-1}),
\end{aligned}
\end{equation}
with $\bar c$ and $\bar c_1$ defined in (\ref{jul620}) and (\ref{20apr104}), respectively.
In particular, this is how the term $\bar c_1$ appears in Theorem~\ref{conj-jun1102-19}.

\subsection*{The nonlinear contribution}

It remains to find the asymptotics of the term~$Y_\ell$ in (\ref{20apr202}), 
and that computation is quite long. The first step is a series
of simplifications.  We start by approximating
$Q(\tau,\eta)$ in~(\ref{20apr202})~by~$\eta$: set
\begin{equation}\label{19jul2642}
\begin{aligned}
\bar Y_\ell=-\frac{1}{\sqrt{4\pi}}
\int_0^\infty\int_0^\infty e^{3\tau/2-\eta\exp(\tau/2)}r_\ell^2(\tau,\eta)\eta 
d\tau d\eta.
\end{aligned}
\end{equation}
\begin{lemma}\label{conj-aug102}
There exists $\gamma>0$ so that
\begin{equation}\label{19aug102}
Y_\ell-\bar Y_\ell=O(\ell^{-\gamma}),\hbox{ as $\ell\to+\infty$}.
\end{equation}
\end{lemma}
This lemma is proved in Section~\ref{sec:proof-of-L52}.
The next approximation involves going back to the original  space-time
variables and restricting the spatial integration to "relatively short" distances 
$x\le (t+1)^\delta$ with $\delta>0$ small.
\begin{lemma}\label{lem-aug102}
There exists $\delta>0$ and $\gamma>0$ so that
\begin{equation}\label{19aug104}
|\bar Y_\ell-\bar Y_\ell^{(1)}|=O(\ell^{-\gamma}),~~\hbox{as $\ell\to+\infty$,}
\end{equation}
where 
\begin{equation}\label{19jun1242}
\begin{aligned}
\bar Y_\ell^{(1)}&=-\frac{1}{\sqrt{4\pi}}
\int_0^\infty\int_0^{(t+1)^\delta} xe^{-x}z_\ell^2(t,x) 
\farc{dx dt}{{(1+t)^{3/2}}},
\end{aligned}
\end{equation}
and $z_\ell(t,x)$ is the solution to
\begin{equation}\label{19jul2646bis}
\pdr{z_\ell}t-\pdrr{z_{\ell}}{x}-\frac{3}{2(t+1)}\Big(z_\ell-\pdr{z_\ell}{x}\Big)
+e^{-x}z_\ell^2=0,
\end{equation}
with $z_\ell(0,x)=\psi_0(x-\ell)$.
\end{lemma}
Note that (\ref{19jul2646bis}) is simply (\ref{jun2064}) with a slightly different 
notation. 
The functions $z_\ell(t,x)$ and $r_\ell(\tau,\eta)$ are related by
\begin{equation}\label{20jan802}
z_\ell(t,x)= 
\sqrt{t+1}r_\ell(\log(t+1),x/\sqrt{t+1}).
\end{equation}
This lemma is proved in Section~\ref{sec:proof-L53} below. It is easy to see why 
(\ref{19aug104}) should be true: we expect that
the function~$z_\ell(t,x)$ behaves roughly as 
\begin{equation}\label{20apr302}
z_\ell(t,x)\sim r_\infty(\ell)x\sim\bar c\ell x,
\hbox{ for $x\gg 1$ and $t\gg 1$.}
\end{equation}
With this input, a back of the envelope computation shows
that the integral in (\ref{19jun1242}) over $x\ge (t+1)^\delta$ is, indeed, small, due to the
exponential decay factor.

The next approximation is to discard the "short times" from the time integration,
looking only at times $t\ge \ell^{2-\alpha}$ with some $\alpha>0$ sufficiently small. 
This is, again, quite expected: as the integration is now over the region
$x\le (t+1)^\delta$, and
$z_\ell(t,x)$ is initially supported at $x=\ell$, it takes the time of roughly $t\sim\ell^2$
to populate the domain of integration, so shorter times can be discarded.  
\begin{lemma}\label{lem-jan1402}
Let $\alpha>0$ sufficiently small. Then, there exists $\gamma>0$ so that
\begin{equation}\label{20jan1402}
|\bar Y_\ell^{(2)}-\bar Y_\ell^{(1)}|=O(\ell^{-\gamma}),~~\hbox{as $\ell\to+\infty$,}
\end{equation}
where 
\begin{equation}\label{20jan1404}
\begin{aligned}
\bar Y_\ell^{(2)}&=-\frac{1}{\sqrt{4\pi}}
\int_{\ell^{2-\alpha}}^\infty\int_0^{(t+1)^\delta} xe^{-x}z_\ell^2(t,x) 
\farc{dx dt}{{(1+t)^{3/2}}},
\end{aligned}
\end{equation}
\end{lemma}
This lemma is proved in Section~\ref{sec:proof-lemma-jan1402} below.

Lemmas~\ref{lem-aug102} and~\ref{lem-jan1402} allow 
us to focus on the integration in the region
$0<x<(t+1)^\delta$ and times~$t>\ell^{2-\alpha}$, 
and consider~$z_\ell(t,x)$  as the solution to (\ref{19jul2646bis})
for $0<x<(t+1)^\delta$, 
with the initial condition $z_{\ell}(0,x)=0$ and the boundary condition at 
$x=(t+1)^\delta$ coming from the 
"outer" solution 
in the self-similar variables:
\begin{equation}\label{19jun1240}
z_{\ell}(t,(1+t)^\delta)=(t+1)^{1/2}r_\ell(\log(t+1),(t+1)^{-1/2+\delta}).
\end{equation}

The analysis of the behavior of $r_\ell(\tau,\eta)$ for $\eta\sim O(1)$ and $\tau\to+\infty$
that we have done so far,
would only give information for~$x\sim\sqrt{t}$, not $x\sim t^\delta$ with~$\delta<1/2$,
as this point corresponds to 
\[
\eta\sim t^{-(1/2-\delta)}=\exp(-(1/2-\delta)\tau)\ll 1.
\] 
To bridge this gap, we need   the following crucial lemma that is proved 
in Section~\ref{sec:hwk}. It refines the informal asymptotics in (\ref{20apr302})
by an extremely important Gaussian factor that interpolates the short time
and the long time behavior. 
\begin{lemma}\label{conj-jul2902}
There exists $\alpha>0$ sufficiently small and 
a constant $K>0$ so that for all $t>\ell^{2-\alpha}$ 
we have
\begin{equation}\label{19jun1240bis}
\farc{z_\ell(t,(t+1)^\delta)}{(t+1)^\delta }\le
(\bar c\ell+K\ell^{1-\delta}) e^{-\ell^2/4(t+1)},
\end{equation}
and for all $\ell^{2-\alpha}<t<\ell^{2+\alpha}$ we have
\begin{equation}\label{20apr802}
\farc{z_\ell(t,(t+1)^\delta)}{(t+1)^\delta }\ge 
(\bar c\ell-K\ell^{1-\delta}) e^{-\ell^2/4(t+1)}.
\end{equation}
\end{lemma}

Lemma~\ref{conj-jul2902}, proved in Section~\ref{sec:hwk}
allows us to look at upper and lower solutions $z_\pm$ for $z_\ell$ in the region 
$x\le (1+t)^\delta$ as the solutions 
to  
\begin{equation}\label{19jun1248}
\pdr{z_{\pm}}t-\pdrr{z_{\pm}}{x}-\frac{3}{2(t+1)}\Big(z_{\pm}-\pdr{z_{\pm}}{x}\Big)+e^{-x}z_{\pm}^2=0,~~x<(t+1)^\delta,
\end{equation}
with the boundary condition
\begin{equation}\label{19jun1244}
z_{\pm}(t,(1+t)^\delta)=
\big(\bar c\ell  \pm  K\ell^{1-\delta}\big) e^{-\ell^2/(4(t+1))}
(1+t)^\delta.
\end{equation}
The initial condition for $z_-(t,x)$ at $t=\ell^{2-\alpha}$ 
is $z_-(\ell^{2-\alpha},x)=0$ for all $x<(1+\ell)^\delta$,
and for $z_+(t,x)$ it is~$z_+(\ell^{2-\alpha},x)=Ce^{-c\ell^{\alpha'}}$ for all $x<(1+\ell)^\delta$,
with some $\alpha'>0$, which
comes from (\ref{20apr210}) below.
We have then the inequality 
\begin{equation}\label{19jul2952}
\bar Y_\ell^+\le \bar Y_\ell^{(2)}\le  \bar Y_\ell^-,
\end{equation}
with
\begin{equation}\label{19jul2954}
\bar Y_\ell^\pm=-\frac{1}{\sqrt{4\pi}}\int_{\ell^{2-\alpha}}^\infty
\int_0^{(t+1)^\delta} xe^{-x}z_{\pm}^2(t,x) 
\farc{dx dt}{(1+t)^{3/2}}. 
\end{equation}
%
%
The last step in the proof of Proposition~\ref{conj-jun1104-19}. 
 is to prove the following.
\begin{lemma}\label{conj-aug104}
There exists 
$\delta>0$ so that 
\begin{equation}\label{aug106}
\bar Y_\ell^\pm=
-\bar c\log\ell-\bar c\log\bar c+k_0\bar c+\frac{\bar c}{2}+O(\ell^{-\delta})
\hbox{ as $\ell\to+\infty$,}
\end{equation}
with the constant $k_0$ as in (\ref{20mar3114}). 
\end{lemma}
The proof of Lemma~\ref{conj-aug104} is presented in Sections~\ref{sec:proof-L55} and \ref{sec:proofofL5.5}. 


Now, putting together (\ref{jan806}) and (\ref{aug106}) gives
\begin{equation}\label{jan810}
r_\infty(\ell)=\bar c\ell+2\bar c \log\ell
+\farc{3}{2}g_\infty\bar c+\bar c_1-\bar c\log\bar c+k_0\bar c+\frac{\bar c}{2}+O(\ell^{-\delta}),
\end{equation}
with some $\delta>0$, finishing the proof of Proposition~\ref{conj-jun1104-19}.

\section{Proofs of auxiliary lemmas}\label{sec:aux}

This section contains the proofs of the auxiliary lemmas used in Section~\ref{sec:5}
to prove Proposition~\ref{conj-jun1104-19}, except for Lemma~\ref{conj-jul2902} that is proved
in Section~\ref{sec:hwk}.

\subsection{The proof of Lemma~\ref{lem-aug102}}\label{sec:proof-L53}

First, we go back to the original, non-self-similar variables: write~
\[
r_\ell(\tau,\eta)=w_\ell(e^\tau-1,\eta e^{\tau/2}),
\]
so that the function $w_\ell(t,x)$ satisfies
\begin{equation}\label{19jul2640}
\begin{aligned}
\pdr{w_\ell}{t}-\pdrr{w_\ell}{x}-\farc{w_\ell}{t+1}+\farc{3}{2(t+1)}
\pdr{w_\ell}{x}+\sqrt{t+1}e^{-x}w_\ell^2=0,
\end{aligned}
\end{equation}
with the initial condition $w_\ell(0,x)=\psi_0(x-\ell)$. 
To revert back even more, we write $w_\ell=z_\ell/\sqrt{t+1}$, with the function
$z_\ell$ that satisfies (\ref{19jul2646bis})
with $z_\ell(0,x)=\psi_0(x-\ell)$. 
Then, we can re-write $\bar Y_\ell$ as
\begin{equation}\label{19jun1236}
\begin{aligned}
\bar Y_\ell&=-
\int_0^\infty\int_0^\infty (t+1)^{3/2} e^{-x}w_\ell^2(t,x)\frac{x}{\sqrt{t+1}}
\farc{dt}{t+1} \farc{dx}{\sqrt{t+1}}=-
\int_0^\infty\int_0^\infty xe^{-x}w_\ell^2(t,x) 
\farc{dx dt}{\sqrt{t+1}}\\
&=-
\int_0^\infty\int_0^\infty xe^{-x}z_\ell^2(t,x) 
\farc{dx dt}{(t+1)^{3/2}}.
\end{aligned}
\end{equation}
Next, we split $\bar Y_\ell$ into two terms, corresponding to small and large $x$:
\begin{equation}\label{19jul2648}
\begin{aligned}
\bar Y_\ell&=-
\int_0^\infty\int_0^{(t+1)^\delta} xe^{-x}z_\ell^2(t,x) 
\farc{dx dt}{(t+1)^{3/2}}-
\int_0^\infty\int_{(t+1)^\delta}^{\infty} xe^{-x}z_\ell^2(t,x) \farc{dx dt}{(t+1)^{3/2}}=
\bar Y_{\ell}^{(1)}+\bar Y_{\ell}^{(c)},
\end{aligned}
\end{equation}
with $\delta\in(0,1/2)$ to be chosen, and split the second term  again:
\begin{equation}\label{19jul2654}
\begin{aligned}
\bar Y_\ell^{(c)}&= 
-\int_0^T\int_{(t+1)^\delta}^{\infty} xe^{-x}z_\ell^2(t,x) \farc{dx dt}{(t+1)^{3/2}}
-\int_T^\infty\int_{(t+1)^\delta}^{\infty} xe^{-x}z_\ell^2(t,x) \farc{dx dt}{(t+1)^{3/2}}=
\bar Y_{\ell}^{(c1)}+\bar Y_{\ell}^{(c2)},
\end{aligned}
\end{equation}
with a large $T$ to chosen. We estimate the term $\bar Y_\ell^{(c2)}$ 
using the simple upper bound
\begin{equation}\label{19jul2650}
z_\ell(t,x)\le C(t+1)^{3/2},
\end{equation}
so that
\begin{equation}\label{19jul2652}
\begin{aligned}
|\bar Y_\ell^{(c2)}|&= 
\int_T^\infty\int_{(t+1)^\delta}^{\infty} xe^{-x}z_\ell^2(t,x) \farc{dx dt}{(t+1)^{3/2}}
\le C\int_T^\infty (t+1)^{5/2}e^{-(t+1)^\delta}dt\le C(1+T)^3\exp(-T^\delta).
\end{aligned}
\end{equation}
To estimate $\bar Y_\ell^{(c1)}$ we slightly refine (\ref{19jul2650}) to 
\begin{equation}\label{19jul2658}
z_\ell(t,x)\le C(t+1)^{3/2}\bar z(t,x-\farc{3}{2}\log(t+1)).
\end{equation}
Here, $\bar z(t,x)$ is the solution to the heat equation on the whole line
\begin{equation}\label{19jul2656}
\pdr{\bar z}{t}=\pdrr{\bar z}{x},~~\bar z(0,x)=\psi_0(x-\ell).
\end{equation}
It follows that  
\begin{equation}\label{19jul2662}
z_\ell(t,x)\le C(t+1)e^{-c\ell},\hbox{ for $0\le t\le \ell$ and $0<x<\ell/2$},
\end{equation}
so that for $0\le t\le \ell$ we have
\begin{equation}\label{19jul2660}
\int_{(t+1)^\delta}^{\infty} xe^{-x}z_\ell^2(t,x)dx\le 
\int_{0}^{\ell/2} xe^{-x}z_\ell^2(t,x)dx+\int_{\ell/2}^\infty xe^{-x}z_\ell^2(t,x)dx
\le C(t+1)^3e^{-c\ell}.
\end{equation}
Thus, if we take $T=\ell$, then 
\begin{equation}\label{19jul2664}
\begin{aligned}
|\bar Y_\ell^{(c1)}|&= 
\int_0^T\int_{(t+1)^\delta}^{\infty} xe^{-x}z_\ell^2(t,x) \farc{dx dt}{(t+1)^{3/2}}
\le C\int_0^T (t+1)^{3}e^{-C\ell}dt\le C(\ell+1)^4\exp(-c\ell).
\end{aligned}
\end{equation}
Using $T=\ell$ also in (\ref{19jul2652}), we obtain
\begin{equation}\label{19jul2668}
\begin{aligned}
|\bar Y_\ell^{(c)}|&\le C(\ell+1)^4\exp(-c\ell^\delta),
\end{aligned}
\end{equation}
finishing the proof of Lemma~\ref{lem-aug102}.

\subsection{The proof of Lemma~\ref{lem-jan1402}}\label{sec:proof-lemma-jan1402}

The difference 
\begin{equation}\label{20jan1406}
\bar Y_\ell^{(1)}-\bar Y_\ell^{(2)}=-\frac{1}{\sqrt{4\pi}}
\int_0^{\ell^{2-\alpha}}\int_0^{(t+1)^\delta} xe^{-x}z_\ell^2(t,x) 
\farc{dx dt}{{(1+t)^{3/2}}}
\end{equation}
can be estimated using (\ref{19jul2658})  and a variant of (\ref{19jul2662}):
\be\label{20apr210}
z_\ell(t,x)\le C(t+1)^{3/2}e^{-c\ell^{\alpha}},~
\hbox{ for $t\le \ell^{2-\alpha}$ and $x\le (t+1)^\delta\le \ell^{(2-\alpha)\delta}$. }
\ee
This gives
\begin{equation}\label{20jan1406bis}
\begin{aligned}
|\bar Y_\ell^{(1)}-\bar Y_\ell^{(2)}|
&\le C\int_0^{\ell^{2-\alpha}}
\int_0^{(t+1)^\delta} xe^{-x}{{(1+t)^{3/2}}}e^{-c\ell^{\alpha}}
{dx dt}
\le C\ell^5 e^{-c\ell^\alpha},
\end{aligned}
\end{equation}
and (\ref{20jan1402}) follows.

\subsection{The proof of Lemma~\ref{conj-aug104}: an informal argument}\label{sec:proof-L55}

Let us explain informally why Lemma~\ref{conj-aug104} is true, before giving the proof
in Section~\ref{sec:proofofL5.5}  below. We expect that
the functions~$z_{\pm}(t,x)$, solutions to (\ref{19jun1248})-(\ref{19jun1244}),
converge   as $t\to+\infty$ to~$\bar z_\pm(x)$, the solutions to 
\begin{equation}\label{19jun1248bis}
-\bar z_\pm''+e^{-x}\bar z_\pm^2=0,
\end{equation}
with the boundary condition
\begin{equation}\label{19jun1250}
\bar z_x(+\infty)=r_\pm:=\bar c\ell\pm K\ell^{1-\delta},
\end{equation}
with the constant $K$ as in Lemma~\ref{conj-jul2902}. 
Note that $\bar z_\pm(x)$ has an explicit form
\begin{equation}\label{19jun1254}
\bar z_\pm(x)=r_\pm \bar z_0(x-\log r_\pm),
\end{equation}
with $\bar z_0$ that solves
\begin{equation}\label{19jun1256}
-\bar z_{0}''+e^{-x}\bar z_0^2=0,
\end{equation}
with the slope at infinity
\begin{equation}\label{19jun1258}
\bar z_{0,x}(+\infty)=1.
\end{equation}
More precisely, we expect that $z(t,x)$ is well approximated by taking
\begin{equation}\label{20jan2310}
r_\pm(t)=r_\pm e^{-\ell^2/(4t)},
\end{equation}
leading to the asymptotics
\begin{equation}\label{19jun1252}
z_{\pm}(t,x)\approx r_\pm(t)\bar z_0(x-\log r_\pm(t))= \bar z_\pm\Big(x+\farc{\ell^2}{4t}\Big)e^{-\ell^2/(4t)},
\end{equation}
that holds for all $t>0$, in the sense that 
\begin{equation}\label{19jul2678}
\begin{aligned}
\bar Y_\ell^{\pm}&=-\farc{1}{\sqrt{4\pi}}
\int_{\ell^{2-\alpha}}^\infty
\int_0^{(t+1)^\delta} xe^{-x}\bar z_\pm^2\Big(x+\frac{\ell^2}{4t}\Big) e^{-\ell^2/(2t)}
\farc{dx dt}{(t+1)^{3/2}}+O(\ell^{-\gamma})\\
&=-\farc{1}{\sqrt{4\pi}\ell}
\int_{\ell^{-\alpha}}^\infty
\int_0^{\infty} xe^{-x}\bar z_\pm^2\Big(x+\frac{1}{4t}\Big) e^{-1/(2t)}
\farc{dx dt}{t^{3/2}}+O(\ell^{-\gamma})\\
&=-\farc{1}{\sqrt{4\pi}\ell}
\int_{0}^\infty
\int_0^{\infty} xe^{-x}\bar z_\pm^2\Big(x+\frac{s}{4}\Big) e^{-s/2}
\farc{dx ds}{\sqrt{s}}+O(\ell^{-\gamma})\\
&=-\farc{1}{\sqrt{4\pi}\ell}
\int_{0}^\infty
\int_{s/4}^{\infty} \Big(x-\farc{s}{4}\Big)e^{-x}\bar z_\pm^2(x) e^{-s/4}
\farc{dx ds}{\sqrt{s}}+O(\ell^{-\gamma})\\
&=-\farc{r_\pm^2}{\sqrt{\pi}\ell}
\int_{0}^\infty
\int_{s}^{\infty}(x-s)e^{-x}\bar z_0^2(x-\log r_\pm) e^{-s}
\farc{dx ds}{\sqrt{s}}+O(\ell^{-\gamma})
\\
&=-\farc{r_\pm}{\sqrt{\pi}\ell}
\int_{0}^\infty
\int_{s-\log r_\pm}^{\infty}(x+\log r_\pm-s)e^{-x}\bar z_0^2(x) e^{-s}
\farc{dx ds}{\sqrt{s}}+O(\ell^{-\gamma})\\
&=-\farc{r_\pm}{\sqrt{\pi}\ell}
\int_{0}^\infty G(\log r_\pm-s) e^{-s}
\farc{ds}{\sqrt{s}}+O(\ell^{-\gamma})
~~~\hbox{ as $\ell\to+\infty$},
\end{aligned}
\end{equation}
with some $\gamma>0$, and
\begin{equation}\label{20apr308}
G(q)=\int_{-q}^\infty (x+q)e^{-x}\bar z_0^2(x)dx. 
\ee
Note that 
\begin{equation}\label{20apr502}
|G(q)|\le C(1+|q|).
\ee
It follows that 
\be\label{20apr504}
\Big|\int_{\log r_\pm/2}^\infty G(\log r_\pm-s) e^{-s}
\farc{ds}{\sqrt{s}}\Big|\le C\int_{\log r_\pm/2}^\infty |\log r_\pm-s| e^{-s}
\farc{ds}{\sqrt{s}}\le \frac{C}{r_\pm^\gamma}\le \farc{C}{\ell^\gamma},
\ee
with $\gamma>0$ sufficiently small. Therefore, we have
\begin{equation}\label{20apr504bis}
\begin{aligned}
\bar Y_\ell^{\pm}&= -\farc{r_\pm}{\sqrt{\pi}\ell}
\int_{0}^{\log r_\pm/2} G(\log r_\pm-s) e^{-s}
\farc{ds}{\sqrt{s}}+O(\ell^{-\gamma})
~~~\hbox{ as $\ell\to+\infty$}.
\end{aligned}
\ee
To analyze the asymptotic behavior of $G(r)$ for $r\gg 1$, note that the solution to (\ref{19jun1256}) with the normalization
(\ref{19jun1258}) 
is given by
\begin{equation}\label{19jul2914}
\bar z(x)=e^{x}U(x),
\end{equation}
where $U(x)$ is the Fisher-KPP minimal speed wave, 
solution to~(\ref{jul616}), with the normalization (\ref{20mar3114}). 
%
Hence, the function $\bar z_0(x)$ has the asymptotics
\begin{equation}\label{19aug108}
\bar z_0(x)=x+k_0+O(e^{-\omega x}),~~\hbox{ as $x\to+\infty$},
\end{equation}
with $k_0$ as in (\ref{20mar3114}), and some $\omega>0$, and
\begin{equation}\label{19jul2916}
\bar z_0(x)\sim e^{x} \hbox{ as $x\to-\infty$}.
\end{equation}
In addition, we have from (\ref{19jun1256})-(\ref{19jun1258}) that
\begin{equation}\label{20jan2112}
\int_{-\infty}^\infty e^{-x}\bar z_0^2(x)dx=1,
\end{equation}
and
\begin{equation}\label{20jan2114}
\begin{aligned}
\int_{-\infty}^\infty x e^{-x}\bar z_0^2(x)dx&
= \int_{-\infty}^\infty x\bar z_0''(x)dx=\lim_{m\to+\infty}\int_{-m}^m x\bar z_0''(x)dx\\
&=
\lim_{m\to+\infty}\big(m\bar z_0'(m)+m\bar z_0'(-m)-\bar z_0(m)+\bar z_0(-m))=-k_0.
\end{aligned}
\end{equation}
Therefore, we can write $G(q)$ as 
\begin{equation}\label{20apr510}
G(q)=\int_{-\infty}^\infty (x+q)e^{-x}\bar z_0^2(x)dx -
\int_{-\infty}^{-q} (x+q)e^{-x}\bar z_0^2(x)dx=q-k_0-\int_{-\infty}^{-q} (x+q)e^{-x}\bar z_0^2(x)dx.
\ee
For the last integral in the right side, we deduce from (\ref{19jul2916}) that 
\begin{equation}\label{19jul2672}
\begin{aligned}
\int_{-\infty}^{-q} (x+q)e^{-x}
\bar z_0^2(x)dx= \int_{-\infty}^{-q} (x+q)e^{x}dx+O(e^{-(1+\delta) q})=-e^{-q}+O(e^{-(1+\delta) q}),~~\hbox{ as $q\to+\infty$}.
\end{aligned}
\end{equation}
It follows that for $s<\log r_\pm/2$ we have
\begin{equation}\label{20apr512}
G(\log r_\pm -s)=\log r_\pm -s-k_0-\farc{1}{r_\pm}e^{s}+O(e^{-(1+\delta)(\log r_\pm-s)})=\log r_\pm -s-k_0+O(\ell^{-\gamma}).
\ee
%
Hence, we have
\begin{equation}\label{20apr504bis2}
\begin{aligned}
\bar Y_\ell^{\pm}&= -\farc{r_\pm}{\sqrt{\pi}\ell}
\int_{0}^{\log r_\pm/2} [\log r_\pm -s-k_0] 
e^{-s}\farc{ds}{\sqrt{s}}+O(\ell^{-\gamma})\\
&=-\farc{r_\pm}{\sqrt{\pi}\ell}
\int_{0}^{\infty} [\log r_\pm -s-k_0] 
e^{-s}\farc{ds}{\sqrt{s}}+O(\ell^{-\gamma})+O(\ell^{-\gamma})\\
&=-\farc{r_\pm}{\sqrt{\pi}\ell}\big([\log r_\pm-k_0]\Gamma(1/2)-\Gamma(3/2)\big)+O(\ell^{-\gamma})
=-\farc{r_\pm}{ \ell}\big(\log r_\pm-k_0 -\frac{1}{2}\big)+O(\ell^{-\gamma})\\
&=-\bar c\log\ell-\bar c\log\bar c+\bar c k_0+\farc{\bar c}{2}+O(\ell^{-\gamma}).
\end{aligned}
\ee
We used above the fact that  $\Gamma(1/2)=\sqrt{\pi}$ and $\Gamma(3/2)=\sqrt{\pi}/2$.

\subsection{The proof of Lemma~\ref{conj-aug104} }\label{sec:proofofL5.5}


Let us now proceed with the actual proof of  Lemma~\ref{conj-aug104}.
The computation in the previous section relied crucially on approximation (\ref{19jun1252}), and the main
step is to justify it. 
The function $z_{+}(t,x)$ satisfies (\ref{19jun1248})-(\ref{19jun1244}) (we drop the sub-script $+$ in this section): 
\begin{equation}\label{19jul2906}
\pdr{z}t-\pdrr{z}{x}-\frac{3}{2(t+1)}\Big(z-\pdr{z}{x}\Big)+e^{-x}z^2=0,~~x<(t+1)^\delta,
\end{equation}
for $t\ge \ell^{2-\alpha}$, with the boundary condition
\begin{equation}\label{19jul2908}
z(t,(t+1)^\delta)=
r_+  e^{-\ell^2/(4(t+1))}(t+1)^\delta,
\end{equation}
with $r_+=\bar c\ell+K\ell^{1-\delta}$, and the initial condition  
\be\label{20apr516}
z(t=\ell^{2-\alpha},x)=e^{-c\ell^{\alpha'}},~~~~x<(t+1)^\delta,
\ee
as in the upper bound (\ref{20apr210}), with some $\alpha'>0$. 
 
Our goal is to find a more explicit super-solution than the solution to (\ref{19jul2906})-(\ref{20apr516}). 
Motivated by (\ref{19jun1254}), let us define 
\begin{equation}\label{jul2910}
\psi(t,x)=e^{\zeta(t)}\bar z_0(x-\zeta(t)),
\end{equation}
with 
\begin{equation}\label{19jul2940}
\zeta(t)=\log r_+ +\log\Big(1+\frac{\log\ell}{\ell^{\delta}}\Big)-\farc{\ell^2}{4(t+1)},
\end{equation}
so that
\be\label{20apr710}
\bal
\psi(t,x)&=r_+\Big(1+\frac{\log\ell}{\ell^{\delta}}\Big) e^{-\ell^2/(4(t+1))}
\bar z_0\Big(x-\log r_\pm-\log\Big(1+\frac{\log\ell}{\ell^{\delta}}\Big)
+\farc{\ell^2}{4(t+1)}\Big)\\
&=\Big(1+\frac{\log\ell}{\ell^{\delta}}\Big) e^{-\ell^2/(4(t+1))}
\bar z_+\Big(x-\log\Big(1+\frac{\log\ell}{\ell^{\delta}}\Big)
+\farc{\ell^2}{4(t+1)}\Big),
\end{aligned}
\ee
with $\bar z_+(x)$, as in (\ref{19jun1254}). 
\begin{lemma}\label{lem-20apr702}
There exists $C>0$ sufficiently large, and $\lambda>0$ sufficiently small, so that 
\begin{equation}\label{20apr604bis}
z(t,x)\le \psi(t,x)+\farc{C}{\ell^{1/4}(t+1)^\lambda},~~\hbox{ for $t>t_{in}=\ell^{2-\alpha}$ 
and $|x|<(t+1)^\delta$.}
\end{equation}
\end{lemma}
\subsection*{The proof of Lemma~\ref{lem-20apr702}}

Our goal will be to show that we can choose $\lambda>0$, $\gamma>0$ and $C>0$ so that
\be\label{20apr604}
s(t,x)=z(t,x)-\psi(t,x)\le 
\bar s(t,x)=\farc{C}{\ell^{1/4}(t+1)^\lambda}\cos\Big(\farc{x}{(t+1)^{\delta+\gamma}}\Big),~~|x|<(t+1)^\delta.
\end{equation}
Let us first show that with the choice of the shift $\zeta(t)$ as in (\ref{19jul2940}), at the boundary $x=(t+1)^\delta$ we have
\begin{equation}\label{19jul2912}
\psi(t,(t+1)^\delta)\ge z(t, (t+1)^\delta),
\end{equation}
so that 
\begin{equation}\label{20apr606}
s(t,(t+1)^\delta)\le 0\le \bar s (t,(t+1)^\delta)\hbox{ for $t\ge T_0$. }
\ee
At this point we have
\begin{equation}\label{19jul2942}
\psi(t,(1+t)^\delta)=r_+\Big(1+\frac{\log\ell}{\ell^\delta}\Big) 
e^{-{\ell^2}/{4(t+1)}}\bar z_0\Big((t+1)^\delta-\log r_+  -\log\Big(1+
\frac{\log\ell}{\ell^{\delta}}\Big) 
+\farc{\ell^2}{4(t+1)}\Big).
\end{equation}
It follows from (\ref{19aug108}) that there exists $k_1$ so that
\begin{equation}\label{19aug110}
\bar z_0(x)\ge x+k_1\hbox{ for all $x\ge 0$}.
\end{equation}
Therefore, if 
\[
\farc{\ell^2}{4(t+1)}\ge \log r_+ +\log\Big(1+\frac{\log\ell}{\ell^{\delta}}\Big)-k_1,
\]
then, as $\bar z_0(x)$ is increasing,  we automatically have from (\ref{19jul2942}) that 
\begin{equation}\label{19jul2944}
\begin{aligned}
\psi(t,(1+t)^\delta)&\ge r_+ 
e^{-{\ell^2}/{4(t+1)}}\Big((t+1)^\delta-\log r_+  
-\log\Big(1+\frac{\log\ell}{\ell^{\delta}}\Big)+\farc{\ell^2}{4(t+1)}+k_1 \Big)\\
& \ge r_+ 
e^{-{\ell^2}/{4(t+1)}}(t+1)^\delta.
\end{aligned}
\end{equation}
On the other hand, if
\begin{equation}\label{19jul3104}
\farc{\ell^2}{4(t+1)}<\log r_++\log\Big(1+\frac{\log\ell}{\ell^{\delta}}\Big)-k_1,
\end{equation}
and $\ell$ is sufficiently large, then we have 
\begin{equation}\label{19aug112}
t>d\farc{\ell^2}{\log\ell}, 
\end{equation}
with a sufficiently small $d>0$,
and $(t+1)^\delta$ dominates the other terms in the argument inside $\bar z_0$ in the right side
of (\ref{19jul2942}). 
In particular, if (\ref{19jul3104}) holds,  since $\bar z_0$ is increasing, we can use (\ref{19aug110}) to obtain
\begin{equation}\label{19jul3106}
\begin{aligned}
\Big(1+\farc{\log\ell}{\ell^\delta}\Big)&\bar z_0\Big(\!(t+1)^\delta-\log r_+  -\log\big(1+
\frac{\log\ell}{\ell^{\delta}}\big) 
+\farc{\ell^2}{4(t+1)}\Big)\ge \Big(1+\farc{\log\ell}{\ell^\delta}\Big)\bar z_0\big((t+1)^\delta-1-\log r_+ \big)\\
&\ge 
 \Big(1+\farc{\log\ell}{\ell^\delta}\Big)\big((t+1)^\delta-1-\log r_+-k_1\big)
 \ge (t+1)^\delta,
\end{aligned}
\end{equation}
because
\begin{equation}\label{19jul3108}
\begin{aligned}
&\farc{\log\ell}{\ell^\delta}\Big((t+1)^\delta
-1-\log r_+-k_1\Big)
-1-\log r_+-k_1
\ge  c\ell^{\delta/2}\log\ell-C\log\ell>0.
\end{aligned}
\end{equation}
We have used (\ref{19aug112}) above. 
Therefore, if we choose $\zeta(t)$ as in (\ref{19jul2940}), then the comparison at the boundary  (\ref{19jul2912}), indeed, holds.

To look at the other boundary point, $x=-(t+1)^\delta$,  
note that the function $\bar z(x)= Ae^x$ is a super-solution to (\ref{19jul2906}) 
for all~$A>1$ sufficiently large,
and at~$x=(t+1)^\delta$ we have, due to Lemma~\ref{conj-jul2902}:
\begin{equation}\label{19aug122}
z(t,(t+1)^\delta)\le C\ell (t+1)^\delta e^{-\ell^2/(4(t+1)}\le C(t+1)^{1/2+\delta}\le Ae^{(t+1)^\delta}=\bar z((t+1)^\delta),
\end{equation}
if $A$ is sufficiently large. It follows that
\[
z(t,x)\le Ae^x,~~x<(t+1)^\delta,
\]
and, in particular, 
we deduce that 
\begin{equation}\label{19aug124}
z(t,-(t+1)^\delta)\le A e^{-(t+1)^\delta}\le\farc{C}{\ell^{1/4}(t+1)^\lambda}\le \bar s(t,-(t+1)^\delta),~~\hbox{for $t\ge \ell^{2-\alpha}$,}
\end{equation}
and thus
\be\label{20apr702}
s(t,-(t+1)^\delta)\le\bar s(t,-(t+1)^\delta),~~t\ge\ell^{2-\alpha}.
\ee
At the initial time $t_{in}=\ell^{2-\alpha}$ we have the comparison
\be\label{20apr704}
s(t_{in},x)\le z(t_{in},x)\le Ce^{-c\ell^{\alpha'}}\le\farc{C}{\ell^{1/4}(1+t_{in})^\lambda}\le \bar s(t_{in},x),~~\hbox{for $x<(t_{in}+1)^\delta$.}  
\ee
%

Having established the comparison at the boundary and the initial time, we  now show that $\bar s(t,x)$ is a super-solution for the equation
for $s(t,x)$. 
The function $\psi(t,x)$ satisfies 
\begin{equation}\label{19jul2926}
\begin{aligned}
\pdr{\psi}t&-\pdrr{\psi}{x}-\frac{3}{2(t+1)}\Big(\psi-\pdr{\psi}{x}\Big)+e^{-x}\psi^2=\dot\zeta(t)e^{\zeta(t)}\bar z_0(x-\zeta(t))-
\dot\zeta (t)e^{\zeta(t)}\pdr{\bar z_0(x-\zeta(t))}{x}\\
&-e^{\zeta(t)}\pdrr{\bar z_0(x-\zeta(t))}{x}
-\frac{3}{2(t+1)}e^{\zeta(t)}\Big(\bar z_0(x-\zeta(t))-\pdr{\bar z_0(x-\zeta(t))}{x}\Big)+
e^{-x}e^{2\zeta(t)}\bar z_0^2(x-\zeta(t))\\
&=\dot\zeta(t)e^{\zeta(t)}\bar z_0(x-\zeta(t))-
\dot\zeta (t)e^{\zeta(t)}\pdr{\bar z_0(x-\zeta(t))}{x}-\frac{3e^{\zeta(t)}}{2(t+1)} \Big(\bar z_0(x-\zeta(t))-\pdr{\bar z_0(x-\zeta(t))}{x}\Big)\\
&=e^{\zeta(t)}\Big(\dot\zeta(t)-\farc{3}{2(t+1)}\Big) \Big(\bar z_0(x-\zeta(t))-\pdr{\bar z_0(x-\zeta(t))}{x}\Big)\\
&\ge 
-e^{\zeta(t)}\Big(\dot\zeta(t)\pdr{\bar z_0(x-\zeta(t))}{x}
+\farc{3}{2(t+1)} \bar z_0(x-\zeta(t))\Big).
\end{aligned}
\end{equation}
We used the fact that the function $\zeta(t)$ is increasing in $t$ and $\bar z_0(x)$ is increasing in $x$ in the last step above.
Note that
\begin{equation}\label{19jul2960}
\begin{aligned}
\dot\zeta(t)e^{\zeta(t)}&=r_+\Big(1+\farc{\log\ell}{\ell^\delta}\Big)\farc{\ell^2}{4(t+1)^2}e^{-\ell^2/(4(t+1))}\le 
\farc{C\ell^3}{(t+1)^{2}}e^{-\ell^2/(4(t+1))} ,
\end{aligned}
\end{equation}
hence
\begin{equation}\label{19jul2962}
\begin{aligned}
\dot\zeta(t)e^{\zeta(t)}\pdr{\bar z_0(x-\zeta(t))}{x}& \le \farc{C\ell^3}{(t+1)^{2}}e^{-\ell^2/(4(t+1))}\le
 \farc{C}{\ell^{1/2}(t+1)^{1/4}}\farc{\ell^{7/2}}{(t+1)^{7/4}}e^{-\ell^2/(4(t+1))}\\
 &\le  \farc{C}{\ell^{1/2}(t+1)^{1/4}}.
\end{aligned}
\end{equation}
The second term in the right side of (\ref{19jul2926}) for $x<(t+1)^\delta$ 
is bounded  by
\begin{equation}\label{19jul3114}
\begin{aligned}
&\frac{3}{2(t+1)}e^{\zeta(t)} \bar z_0(x-\zeta(t))\le \farc{C\ell}{t+1}e^{-\ell^2/(4(t+1))}\Big( (t+1)^\delta 
+\farc{\ell^2}{t+1}\Big)\\
&=C\Big(\farc{\ell}{(t+1)^{1-\delta}}
+\farc{\ell^3}{(t+1)^2}\Big)e^{-\ell^2/(4(t+1))}
\\
&\le C\Big(\farc{1}{\ell^{1/4}(t+1)^{3/8-\delta}}\cdot\farc{\ell^{5/4}}{(t+1)^{5/8}}
+\farc{1}{\ell^{1/2}(t+1)^{1/4}}\cdot\farc{\ell^{7/2}}{(t+1)^{7/4}}\Big)
e^{-\ell^2/(4(t+1))}\le \farc{C}{\ell^{1/4}(t+1)^{1/4}}.
\end{aligned}
\end{equation}
%
Therefore, the function $s(t,x)=z(t,x)-\psi(t,x)$ satisfies
\begin{equation}\label{19jul3110}
\begin{aligned}
\pdr{s}t&-\pdrr{s}{x}-\frac{3}{2(t+1)}\Big(s-\pdr{s}{x}\Big)+e^{-x}(\psi(t,x)+z(t,x))s\le \farc{C}{\ell^{1/4}(t+1)^{1/4}}
,~~x<(t+1)^\delta.
\end{aligned}
\end{equation}
On the other hand, recall from~\cite{NRR1} that for $\gamma>0$, the function $\bar s(t,x)$ 
satisfies  
\begin{equation}\label{19jul3120}
\pdr{\bar s(t,x)}{t}=-\farc{\lambda}{t+1}\bar s(t,x)+g(t,x),
\end{equation}
with $g(t,x)$ such that
\begin{equation}\label{19jul3122}
|g(t,x)|\le \farc{C(\delta+\gamma)|x|}{\ell^{1/4}(t+1)^{\lambda+\delta+\gamma+1}}\le  \farc{C(\delta+\gamma)\bar s(t,x)}{ (t+1)^{\gamma+1}},
~~\hbox{  for $|x|\le (t+1)^\delta$.}
\end{equation}
In addition, we have 
\begin{equation}\label{19jul3124}
\pdrr{\bar s(t,x)}{x}=-\farc{1}{(t+1)^{2\delta+2\gamma}}\bar s(t,x),
\end{equation}
and 
\begin{equation}\label{19jul3126}
\farc{3}{2(t+1)}\Big(\pdr{\bar s(t,x)}{x}-\bar s(t,x)\Big)\ge -\farc{C}{t+1}\bar s(t,x),
~~\hbox{  for $|x|\le (t+1)^\delta$.}
\end{equation}
We conclude that  
\begin{equation}\label{19jul3128}
\begin{aligned}
\pdr{\bar s}t&-\pdrr{\bar s}{x}-\frac{3}{2(t+1)}\Big(\bar s-\pdr{\bar s}{x}\Big)
+e^{-x}(\psi(t,x)+z(t,x))\bar s\\
&\ge 
-\farc{\lambda}{1+t}\bar s(t,x)+g(t,x)+\farc{1}{(1+t)^{2\delta+2\gamma}}\bar s(t,x)
-\farc{C\bar s(t,x)}{(t+1)}
\ge \farc{C}{(1+t)^{2\delta+2\gamma}}\bar s(t,x)\\\
&\ge \farc{C}{\ell^{1/4}(1+t)^{\lambda +2\delta+2\gamma}}, 
~~\hbox{ for $t\ge \ell^{2-\alpha}$ and $|x|<(t+1)^\delta$,}
\end{aligned}
\end{equation}
provided that $\delta>0$ and $\gamma>0$ are sufficiently small. 
%
Now, putting together the boundary comparisons (\ref{20apr606}) and (\ref{20apr702}),
the initial time comparison (\ref{20apr704}), as well as 
(\ref{19jul3110}) and (\ref{19jul3128}), 
we deduce from the comparison principle that
if $\lambda>0$, $\delta>0$ and $\gamma>0$ are all sufficiently small, then
\begin{equation}\label{19aug136}
s(t,x)\le \bar s(t,x)\hbox{ for $t>t_{in}=\ell^{2-\alpha}$ 
and $|x|<(t+1)^\delta$.}
\end{equation}
This implies (\ref{20apr604}), and (\ref{20apr604bis}) follows, finishing the proof
of Lemma~\ref{lem-20apr702}.~$\Box$


\subsubsection*{The end of the proof of Lemma~\ref{conj-aug104}}

Using  (\ref{20apr604bis}) in (\ref{19jul2954}) gives 
\begin{equation}\label{19jul3132}
\begin{aligned}
-\bar Y_\ell^{(+)}&=\farc{1}{\sqrt{4\pi}}
\int_{\ell^{2-\alpha}}^\infty\int_0^{(t+1)^\delta} xe^{-x}z_\ell^2(t,x) 
\farc{dx dt}{{(t+1)^{3/2}}}\\
&\le \farc{1}{\sqrt{4\pi}}\int_{\ell^{2-\alpha}}^\infty
\int_0^{(t+1)^\delta} xe^{-x}\Big(\psi(t,x) +\farc{C}{\ell^{1/4}(t+1)^\lambda}\Big)^2
\farc{dx dt}{{(t+1)^{3/2}}} 
=I+II+III, 
\end{aligned}
\end{equation}
with the three terms coming from the expansion
\[
\Big(\psi(t,x) +\farc{C}{\ell^{1/4}(t+1)^\lambda}\Big)^2=\psi^2(t,x)+\farc{2C}{\ell^{1/4}(t+1)^\lambda}\psi(t,x)+
\farc{C^2}{\ell^{1/2}(t+1)^{2\lambda}}.
\]
We have, clearly,
\begin{equation}\label{19jul3136}
III\le \farc{C}{\ell^{1/2}}.
\end{equation}
For the second term, we recall that for $x>0$ we have
\begin{equation}\label{19jul3140}
\begin{aligned}
\psi(t,x)&=e^{\zeta(t)}\bar z_0(x-\zeta(t))\le C\ell e^{-\ell^2/(4(t+1))}\bar z_0\Big(x
+\farc{\ell^2}{4(t+1)}\Big)
\\
&\le C \ell e^{-\ell^2/(4(t+1))}\Big(x+ 1+
\farc{\ell^2}{4(t+1)}\Big),
\end{aligned}
\end{equation}
and use this to write
\begin{equation}\label{19jul3138}
\begin{aligned}
II&=\farc{C}{\ell^{1/4}}
\int_{\ell^{2-\alpha}}^\infty\int_0^{(t+1)^\delta} xe^{-x}\psi(t,x)  
\farc{dx dt}{{(t+1)^{3/2+\lambda}}}\\
&\le
\farc{C}{\ell^{1/4}}
\int_0^\infty\int_0^{(t+1)^\delta} xe^{-x} \ell e^{-\ell^2/(4(t+1))}\Big(x+1 
+\farc{\ell^2}{4(t+1)}\Big)
\farc{dx dt}{{(t+1)^{3/2+\lambda}}}\\
&\le \farc{C}{\ell^{1/4}}
\int_1^\infty  \ell e^{-\ell^2/(4t)}\Big(1 
+\farc{\ell^2}{4t}\Big)
\farc{dx dt}{{t^{3/2+\lambda}}}\le \farc{C}{\ell^{1/4}}
\int_0^\infty  \ell e^{-1/(4t)}\Big(1 
+\farc{1}{4t}\Big)
\farc{\ell^2 dx dt}{{(\ell^2t)^{3/2+\lambda}}}\le\farc{C}{\ell^{1/4}}.
\end{aligned}
\end{equation}
For the main term, we need to compute more precisely, and we use expression
(\ref{20apr710}) for $\psi(t,x)$:
\begin{equation}\label{19jul3142}
\begin{aligned}
&I= \farc{1}{\sqrt{4\pi}}\int_{\ell^{2-\alpha}}^\infty
\int_0^{(t+1)^\delta} xe^{-x}\psi^2(t,x)  \farc{dx dt}{{(1+t)^{3/2}}}\\
&=\Big(1+\farc{\log\ell}{\ell^\delta}\Big)^2
\!\int_{\ell^{2-\alpha}}^\infty\int_0^{(t+1)^\delta}\!\!\!\!\!\! \!\!xe^{-x}e^{-\ell^2/(2(t+1))}
\bar z_+^2\Big(x-\log\Big(1+\farc{\log\ell}{\ell^\delta}\Big)+\farc{\ell^2}{4(t+1)}\Big)  \farc{dx dt}{{(t+1)^{3/2}}}\\
&= 
\!\int_{\ell^{2-\alpha}}^\infty\int_0^{(t+1)^\delta}\!\!\!\!\!\! \!\!xe^{-x}e^{-\ell^2/(2(t+1))}
\bar z_+^2\Big(x+\farc{\ell^2}{4(t+1)}\Big)  \farc{dx dt}{{(t+1)^{3/2}}}+O(\ell^{-\delta/2}).
\end{aligned}
\end{equation}
This is simply expression (\ref{19jul2678}) that we have computed in
Section~\ref{sec:proof-L55}, leading to  (\ref{20apr504bis2}):
\begin{equation}\label{20jan2416}
-\bar Y_{\ell}^+\le \bar c\log\ell+\bar c\log\bar c-\bar ck_0
-\farc{\bar c}{2}+O(\ell^{-\gamma})
\end{equation}
This finishes the proof of the upper bound in Lemma~\ref{conj-aug104}.

Proceeding as in the proof of the upper bound, with some minor modifications, 
we can obtain a matching lower bound
\begin{equation}\label{19aug174}
-\bar Y_{\ell}^-\ge \bar c\log\ell+\bar c\log\bar c-\bar c k_0
-\farc{\bar c}{2}+O(\ell^{-\gamma})
\end{equation}
which finishes the proof of Lemma~\ref{conj-aug104}. 

\subsection{The proof of Lemma~\ref{conj-aug102}}\label{sec:proof-of-L52}

We need to show that 
\begin{equation}\label{19aug204}
\begin{aligned}
\tilde Y_\ell&=- 
\int_0^\infty\int_0^\infty e^{3\tau/2-\eta\exp(\tau/2)}r_\ell^2(\tau,\eta)\bar\psi(\eta)e^{-\tau/2}
d\tau d\eta=O(\ell^{-\gamma}),\hbox{ as $\ell\to+\infty$}.
\end{aligned}
\end{equation}
As in (\ref{19jun1236}), we re-write this in terms of $z_\ell(t,x)$ as 
\begin{equation}\label{19aug206}
\begin{aligned}
\tilde Y_\ell&=-
\int_0^\infty\int_0^\infty (t+1)  e^{-x}w_\ell^2(t,x)\bar\psi\Big(\frac{x}{\sqrt{t+1}}\Big)
\farc{dt}{t+1} \farc{dx}{\sqrt{t+1}}\\
&
=-
\int_0^\infty\int_0^\infty e^{-x}z_\ell^2(t,x) \bar\psi\Big(\frac{x}{\sqrt{t+1}}\Big)
\farc{dx dt}{(t+1)^{3/2}}.
\end{aligned}
\end{equation}
We sketch the argument that can be made precise using Lemma~\ref{lem-20apr702} in a straightforward way. 
Let us use approximation (\ref{19jun1252}) again:
\begin{equation}\label{19aug208}
z_{\ell}(t,x)\approx \bar z_\pm(x)e^{-\ell^2/(4(t+1))},
\end{equation}
and insert this into (\ref{19aug206}). This would give
\begin{equation}\label{19aug210}
\begin{aligned}
\tilde Y_\ell&\approx - \int_0^\infty\int_0^\infty e^{-x}\bar z_\pm^2(x)e^{-\ell^2/(2(t+1))} \bar\psi\Big(\frac{x}{\sqrt{t+1}}\Big)
\farc{dx dt}{(t+1)^{3/2}}\\
&  \approx -r_\pm^2 \int_0^\infty\int_0^\infty e^{-x}\bar z_0^2(x-\log r_\pm)e^{-\ell^2/(2(t+1))} \bar\psi\Big(\frac{x}{\sqrt{t+1}}\Big)
\farc{dx dt}{(t+1)^{3/2}}.
\end{aligned}
\end{equation}
Recall that $\bar\psi(\eta)\le C\eta$, so we can roughly estimate
\begin{equation}\label{19aug214}
\begin{aligned}
|\tilde Y_\ell|&  \le Cr_\pm^2 \int_0^\infty\int_{-\infty}^{\infty} 
e^{-x}\bar z_0^2(x-\log r_\pm)e^{-\ell^2/(2(t+1))}\frac{|x|}{\sqrt{t+1}} 
\farc{dx dt}{(t+1)^{3/2}} \\
&\le Cr_\pm \int_0^\infty\int_{-\infty}^{\infty} 
e^{-x}\bar z_0^2(x)e^{-\ell^2/(2(t+1))}\frac{|x|+\log r_\pm}{\sqrt{t+1}} 
\farc{dx dt}{(t+1)^{3/2}}\\
&\le C r_\pm(1+\log r_\pm)\int_0^\infty 
 e^{-\ell^2/(2(t+1))} 
\farc{dt}{(t+1)^{2}}=C r_\pm(1+\log r_\pm)\farc{\ell^2}{\ell^4}\le\farc{C\log\ell}{\ell}.
\end{aligned}
\end{equation}
This finishes the proof of Lemma~\ref{conj-aug102}.~$\Box$


\subsection{The proof of Lemma~\ref{conj-jul2902}}\label{sec:hwk}

We now turn to the proof of Lemma~\ref{conj-jul2902}. 
Let us go back to (\ref{19jul2646bis}):
\begin{equation}\label{20apr806}
\pdr{z_l}{t}-\pdrr{z_l}{x}-\frac3{2(t+1)}\Big(z_l-\pdr{z_l}{x}\Big)+e^{-x}z_l^2=0,
\ee
and undo yet another change of variables: 
\begin{equation}
\label{e1.4}
\bal
v(t,x)=&(t+1)^{-3/2}z_l(t,x+\farc{3}{2}\log (t+1)).
\enbal
\end{equation}
The function $v(t,x)$ satisfies simply
\begin{equation}
\label{20apr808}
\bal
&\pdr{v}t-\pdrr{v}{x}+e^{-x}v^2=0\\
& v(0,x)=\psi_0(x-\ell).
\enbal
\ee
Note that
\begin{equation}\label{20apr810}
\farc{z_\ell(t,(t+1)^\delta)}{(t+1)^\delta}=
\farc{(t+1)^{3/2}v(t,(t+1)^\delta-(3/2)\log(t+1))}{(t+1)^\delta},
\ee
so that our point of interest for $v(t,x)$ is
\be\label{20apr816}
x_\ell(t)=(t+1)^\delta-\farc{3}{2}\log(t+1).
\ee
%

We perform the standard parabolic scaling
\begin{equation}\label{20apr812}
p(s,\xi)=v(\ell^2s,\ell\xi), 
\ee
so that the function $p(s,\xi)$ solves
\begin{equation}\label{e1.6}
\bal
&\pdr{p}s-\pdrr{p}{\xi}+\ell^2e^{-\ell\xi}p^2=0\\
& p(0,\xi)=\psi_0(\ell(\xi-1)).
\enbal
\end{equation}
In the new variables, the point of interest is
\be\label{20apr814}
\xi_\ell(s)=\farc{x_\ell(\ell^2s)}{\ell}=
\frac{(\ell^2s+1)^\delta-(3/2)\log(\ell^2s+1)}{\ell}
\approx\frac{s^\delta}{\ell^{1-2\delta}},
\ee
and the important time scales are $\ell^{-\alpha}\ll s\ll \ell^{\alpha}$, corresponding
to $\ell^{2-\alpha}\ll t\ll \ell^{2+\alpha}$.  In particular, we have
\be\label{20apr1314}
\xi_\ell(s)= 
(1+O(\ell^{-\gamma}))\frac{s^\delta}{\ell^{1-2\delta}},~~
\hbox{ for $\ell^{-\alpha}\ll s\ll \ell^{\alpha}$}. 
\ee

Let us start with the proof of the lower bound (\ref{20apr802}). 
The maximum principle implies that
\[
p(s,\xi)\le C_0=\|\psi_0\|_{L^\infty}.
\]
Therefore, if we take $0<\gamma\ll\delta$ and impose the Dirichlet boundary condition
at $\xi=\ell^{-(1-\gamma)} $, then 
we have 
\be\label{20apr1320}
p(s,\xi)\geq q(s,\xi),
\ee
where $q(s,\xi)$ is the solution to
\begin{equation}\label{e1.8}
\bal
&\pdr{q}{s}-\pdrr{q}{\xi}+C_0\ell^2e^{-\ell^\gamma}q=0,~~s>0,~~\xi>\ell^{-(1-\gamma)},\\
&q(0,\xi)=\psi_0(l(\xi-1)),\\
&q(0,\ell^{-(1-\gamma)})=0.
\enbal
\end{equation}
We see from (\ref{20apr814}) that, as $\gamma<\delta$, we have  
\be\label{20apr1302}
\xi_\ell(s)\ge\farc{s^\delta}{2\ell^{1-2\delta}}\ge 
\frac{\ell^{-\alpha\delta}}{2\ell^{1-2\delta}}
 >\ell^{-(1-\gamma)},~~\hbox{ for $s\ge \ell^{-\alpha}$,}
\ee
provided that $\gamma$ is sufficiently small, 
so that we have not lost our point of interest in this approximation.
An explicit formula for the solution to (\ref{e1.8}) gives
\be\label{20apr820}
q(s,\xi_\ell(s))=\frac{e^{-C_0\ell^2e^{-\ell^\gamma}s}}{\sqrt{4\pi s}}
\int_0^{+\infty}\Big(e^{-(\tilde\xi_\ell(s)-\zeta)^2/(4s)}-
e^{-(\tilde\xi_\ell(s)+\zeta)^2/(4s)}\Big)\psi_0(\ell(\zeta-1))d\zeta,
~~\hbox{ for $s\ge \ell^{-\alpha}$,}
\ee
with
\be\label{20apr1310}
\tilde\xi_\ell(s)=\xi_\ell(s)-\ell^{-(1-\gamma)}\approx\farc{s^\delta}{\ell^{1-2\delta}}
-\farc{1}{\ell^{1-\gamma}}\approx \farc{s^\delta}{\ell^{1-2\delta}}, 
~~\hbox{ for $s\ge \ell^{-\alpha}$.}
\ee
%
%
The restriction on $s$ in (\ref{20apr820}) comes from the requirement that (\ref{20apr1302})
holds, so that $\xi_\ell(s)$ is in the region where $q(s,\xi)$ is defined. 

The integrand in (\ref{20apr820}) is very small for $\|\zeta-1\|\gg \ell^{-1}$ because
of the exponential decay of the initial condition $\psi_0(x)$.  
With this in mind, we write the difference of the exponentials in (\ref{20apr820}) as
\be\label{20apr822}
\bal
e^{-(\tilde\xi_\ell(s)-\zeta)^2/(4s)}-
e^{-(\tilde\xi_\ell(s)+\zeta)^2/(4s)}=
e^{-(\tilde\xi_\ell(s)-\zeta)^2/(4s)}(1-e^{-\tilde\xi_\ell(s)\zeta/s}).
\enbal
\ee
As $\zeta\approx 1$, we have 
\be\label{20apr824}
\farc{\tilde\xi_\ell(s)\zeta}{s}\approx\farc{1}{\ell^{1-2\delta}s^{1-\delta}}\ll 1
~~~\hbox{~~for $s\gg \ell^{(1-2\delta)/(1-\delta)}$.}
\ee
Using the inequality 
\be\label{20apr1304}
1-e^{-u}\geq u-u^2
\ee
for sufficiently small $u$ we have then 
\be\label{20apr826}
q(s,\xi_\ell(s))\ge \frac{e^{-sC_0l^2e^{-l^\gamma}}\tilde\xi_\ell(s)}
{\sqrt{4\pi}s^{3/2}}(1+o(\ell^{-1}))\int_0^{+\infty}
e^{-(\tilde\xi_\ell(s)-\zeta)^2/(4s)}\zeta\Big(1-\frac{\zeta\tilde\xi_\ell(s)}{s}\Big)
\psi_0(\ell (\zeta-1))d\zeta .
\ee
The $o(\ell^{-1})$ correction in the pre-factor comes from $\zeta$ that violate
(\ref{20apr824}), so that we can not use~(\ref{20apr1304}). 
Their contribution   is extremely small since $\psi_0(x)$ is 
decaying exponentially as $x\to-\infty$ and has compact support for~$x>0$. 
Using the straightforward approximations for  the Gaussian and the
factor inside the parenthesis inside the integral in (\ref{20apr826}), 
as well as for the exponential factor in front of the integral, we obtain
\be\label{20apr828}
\bal
q(s,\xi_\ell(s))&\ge \frac{e^{-sC_0l^2e^{-l^\gamma}}\tilde\xi_\ell(s)}
{\sqrt{4\pi}s^{3/2}}(1+O(\ell^{-\gamma}))\int_0^{+\infty}
e^{-1/(4s)}\zeta \psi_0(\ell (\zeta-1))d\zeta\\
&\ge \frac{\tilde\xi_\ell(s)e^{-1/(4s)}}
{\sqrt{4\pi}s^{3/2}}(1+O(\ell^{-\gamma}))\int_0^{+\infty}
 \zeta \psi_0(\ell (\zeta-1))d\zeta\\
 &=
\frac{\tilde\xi_\ell(s)e^{-1/(4s)}}
{\sqrt{4\pi}s^{3/2}\ell}(1+O(\ell^{-\gamma}))\int_{-\ell}^{+\infty}
\psi_0(\zeta)d\zeta=
\frac{\bar c\tilde\xi_\ell(s)e^{-1/(4s)}}
{ s^{3/2}\ell}(1+O(\ell^{-\gamma})).
\enbal
\ee
Going back to (\ref{20apr1314}) and (\ref{20apr1310}), we note that
\be\label{20apr1312}
\tilde\xi_\ell(s)=\xi_\ell(s)-\ell^{-(1-\gamma)}=
(1+O(\ell^{-\gamma}))\farc{s^\delta}{\ell^{1-2\delta}}
-\farc{1}{\ell^{1-\gamma}}\approx \farc{s^\delta}{\ell^{1-2\delta}}= 
(1+O(\ell^{-\gamma}))\xi_\ell(s),
~~\hbox{ for $s\ge \ell^{-\alpha}$.}
\ee
Using this in (\ref{20apr828}) gives
\be\label{20apr1316}
\bal
q(s,\xi_\ell(s))&\ge
\frac{\bar c\xi_\ell(s)e^{-1/(4s)}}
{ s^{3/2}\ell}(1+O(\ell^{-\gamma})),~~~~\hbox{ for $s\ge \ell^{-\alpha}$.}
\enbal
\ee
Unrolling the changes of variables (\ref{20apr810}), (\ref{20apr812}), and (\ref{20apr814}),
and using the bound (\ref{20apr1320}), as well as the approximation (\ref{20apr1314}),
we can re-write (\ref{20apr1316}) as 
\be\label{20apr1318}
\bal
\farc{z_\ell(t,(t+1)^\delta)}{(t+1)^\delta}&=
\farc{(t+1)^{3/2}v(t,(t+1)^\delta-(3/2)\log(t+1))}{(t+1)^\delta}\\
&=\farc{(t+1)^{3/2}p(\ell^{-2}t,\ell^{-1}(t+1)^\delta-(3/2)\ell^{-1}\log(t+1))}{(t+1)^\delta}
=\farc{(t+1)^{3/2}p(\ell^{-2}t,\xi_\ell(\ell^{-2}t)}{(t+1)^\delta}\\
&\ge \farc{(t+1)^{3/2}q(\ell^{-2}t,\xi_\ell(\ell^{-2}t)}{(t+1)^\delta}
\ge \farc{(t+1)^{3/2}\bar c\xi_\ell(\ell^{-2}t)e^{-\ell^2/(4t)}\ell^3}
{(t+1)^\delta(t+1)^{3/2}\ell}(1+O(\ell^{-\gamma}))\\
&=\farc{\bar c\ell^{-2\delta}t^\delta e^{-\ell^2/(4t)}\ell^2}
{(t+1)^\delta \ell^{1-2\delta}}(1+O(\ell^{-\gamma}))=\bar c\ell(1+O(\ell^{-\gamma}))
e^{-\ell^2/(4t)},\hbox{ for $\ell^{2-\alpha}\le t\le \ell^{2+\alpha}$},
\enbal
\ee
which is the lower bound (\ref{20apr802}).

For the upper bound (\ref{19jun1240bis}), 
we again look at the solution $p(s,\xi)$ of \eqref{e1.6}. A simple upper bound for
$p(s,\xi)$ is by the solution to the heat equation on the whole real line:
\begin{equation}\label{20apr1322}
\bal
&\pdr{\bar p}s-\pdrr{\bar p}{\xi}=0,~~x\in\Rm,\\
& \bar p(0,\xi)=\psi_0(\ell(\xi-1)).
\enbal
\end{equation}
Accordingly, we set  
\begin{equation}
\label{e1.11}
\overline\psi_0(\xi)=\bar p(s_\alpha,\xi)=
\frac1{\sqrt{4\pi s_\alpha}}\int_{\R}e^{-(\xi-\zeta)^2/4s_\alpha}\psi_0(\ell(\zeta-1))d\zeta,~~s_\alpha=\ell^{-\alpha}, 
\end{equation}
with a sufficiently small $\alpha$ to be chosen, depending on $\delta$. 
We also have the upper bound
\be\label{20apr1324}
p(s,\xi)\le e^{\ell\xi},
\ee
that follows immediately from the maximum principle. In particular, we have
\be\label{20apr1326}
p(s,-\ell^{-(1-\gamma)})\le e^{-\ell^\gamma},
\ee
whence
\be\label{20apr1328}
p(s,\xi)\leq e^{-\ell^\gamma}+\bar p_1(s,\xi).
\ee
Here, $\bar p_1(s,\xi)$ is the solution to and $\overline v(\tau,\xi)$ solves
\be\label{20apr1330}
\bal 
&\pdr{\bar p_1}{s}-\pdrr{\bar p_1}{\xi}=0,~~s>s_\alpha,~\xi>-\ell^{-(1-\gamma)},\\
&\bar p_1(s_\alpha,\xi)=\overline\psi_0(\xi),\\
&\bar p_1(s,-\ell^{-(1-\gamma)})=0.
\enbal
\ee
Let us note the following properties of the initial condition $\overline \psi_0$. First,
(\ref{e1.11}) implies that it is localized near $\xi=1$, in the sense that
\begin{equation}
\label{e1.12}
0<\overline\psi_0(\xi)\leq C\ell^{-\alpha/2}\ell^{-1}\exp\big\{-\ell^{\alpha-2\beta}\big\},
~~\hbox{ for $|\xi-1|\ge \ell^{-\beta}$.}
\end{equation}
Hence, as soon as $\xi$ departs from a very small neighborhood of 1, 
$\overline\psi_0(\xi)$ is exponentially small in $\ell$. 
Furthermore, the mass of $\overline\psi_0(\xi)$  is
\be\label{20apr1332}
\di\int_\R\overline\psi_0(\xi)d\xi=\di\frac{\bar c}{\ell},
\ee
and its first moment is
\be\label{20apr1334}
\di\int_\R\xi \overline\psi_0(\xi)d\xi=\di\frac{\bar c}{\ell}(1+O(\frac1\ell)),
\ee
because of \eqref{e1.12}.
 
It follows that the function $\bar p_1(s,\xi)$ can  
then be estimated along the same lines as $q(s,\xi)$ in the proof of the lower bound.
This eventually leads to (\ref{19jun1240bis}). $\Box$

\begin{appendix}
\section{Proof of an auxiliary lemma}\label{sec:appendix}

Here we prove an elementary result used in the proof of Lemma~\ref{lem:20_06_1}. 
\begin{lemma}\label{lem-20june1802}
Let $\psi\in \cC^+_{bc}$, $\psi(\cdot)\leq 1$, and 
$g$ be a continuous function such that $g(x)=0$ for $x\leq 0$, $g(x)=1$ for 
$x\geq 1$ and $g(x)=x$ for $x\in [0,1]$.  Define 
$\psi_n(x) = g(x+n)\psi(x)\in \cC^+_c$, then $\hat s[\psi_n]\rightarrow \hat s[\psi]$ as~$n\rightarrow \infty$.
\end{lemma}
{\bf Proof.} As $\psi_n(x)\le \psi(x)$, the comparison principle implies immediately that
\be\label{20jun1802}
\hat s[\psi]\le \hat s[\psi_n],
\ee
and we only need to verify an opposite bound.
Let $u(t,x)$, $u_n(t,x)$ and $\tilde u_n(t,x)$  be the solutions 
to~(\ref{20apr1402}) with the initial conditions 
\be\label{20jun1804}
u(0,x)=\psi(x),~~~u_n(t,x)=\psi_n(x),~~~
\tilde u_n(0,x)=\tilde\psi_n(x):=\psi(x)-\psi_n(x).
\ee
The function $v_n=u_n+\tilde u_n$ satisfies
\be\label{20jun1806}
\pdr{v_n}{t}=\pdrr{u_n}{x}+u_n-u_n^2+\pdrr{\tilde u_n}{x}+\tilde u_n-\tilde u_n^2
\ge \pdrr{v_n}{x}+v_n-v_n^2,
\ee
with the initial condition $v_n(0,x)=u(0,x)$. It follows from the maximum
principle that 
\[
\hbox{$v_n(t,x)\ge u(t,x)$ for all $t\ge 0$ and $x\in\Rm$},
\]
and, in particular,
\be\label{20jun1808}
u(t,x+m(t))\le u_n(t,x+m(t))+\tilde u_n(t,x+m(t)).
\ee
Passing to the limit $t\to+\infty$, we obtain
\be\label{20jun1810}
U(x+\hat s[\psi])\le U(x+\hat s[\psi_n])+U(x+\hat s[\tilde\psi_n]),
\ee
for each $n\in\Nm$ fixed and all $x\in\Rm$. 
Dividing by $x\exp(-x)$ and passing to the limit $x\to+\infty$, keeping
$n\in\Nm$ fixed, gives 
\be\label{20jun1812}
e^{-\hat s[\psi]}\le e^{-\hat s[\psi_n]}+e^{-\hat s[\tilde\psi_n]}.
\ee
As $\tilde\psi_n(x)=0$ for all $x\ge -n+1$, we know that $\hat s[\tilde\psi_n]\to+\infty$
as $n\to+\infty$. Using this and passing to the limit $n\to+\infty$ in (\ref{20jun1812})
leads to
\be\label{20jun1814}
\hat s[\psi]\ge \limsup_{n\to+\infty}\hat s[\psi_n].
\ee
This, together with (\ref{20jun1802}) finishes the proof.~$\Box$

\end{appendix}

\end{document}